\documentclass[EJP, preprint]{ejpecp} 

\usepackage[T1]{fontenc}
\usepackage[utf8]{inputenc}

\usepackage{enumerate} 
\newenvironment{myenumerate}{%
\renewcommand{\theenumi}{(\roman{enumi})}%
\renewcommand{\labelenumi}{\theenumi}%
\begin{list}{\labelenumi}
	{%
	\setlength{\itemsep}{0.4em}%
	\setlength{\topsep}{0.5em}%
	\setlength\leftmargin{2.45em}%
	\setlength\labelwidth{2.05em}%
	\setlength{\labelsep}{0.4em}%
	\usecounter{enumi}%
	}%
	}%
{\end{list}
}
\renewenvironment{enumerate}{
\begin{myenumerate}}%
{\end{myenumerate}}

\usepackage[dvipsnames]{xcolor}
\usepackage{graphicx}
\usepackage{bm}
\usepackage{subcaption}
\usepackage{stmaryrd}

\usepackage{tikz}
\usetikzlibrary{decorations.pathreplacing}
\newcommand{\tikznode}[3][inner sep=0pt]{\tikz[remember
picture,baseline=(#2.base)]{\node(#2)[#1]{$#3$};}}

\usepackage{forest}
\forestset{
  default preamble={
   for tree={circle, draw, 
            minimum size=2.1em, 
            inner sep=1pt} 
  }
}

\usetikzlibrary{positioning}

\usepackage{ytableau}

\SHORTTITLE{$\lambda$-shaped random matrices, $\lambda$-plane trees, and $\lambda$-Dyck paths}

\TITLE{\vspace{-4cm}
$\lambda$-shaped random matrices, $\lambda$-plane trees, and $\lambda$-Dyck paths\support{
This project received financial support from the mobility programme Erasmus+, from the International Office of TU Wien, from the Italian National Group of Mathematical Physics (GNFM-INdAM), from Regione Puglia through the project `Research for Innovation' - UNIBA024,  from Istituto Nazionale di Fisica Nucleare (INFN) through the project `QUANTUM', and from PRIN 2022 project 2022TEB52W-PE1 - `The charm of integrability: from nonlinear waves to random matrices'.}\
    } 

\AUTHORS{%
  Elia~Bisi\footnote{Technische Universit\"at Wien, Austria and Universit\`a degli Studi di Firenze, Italy.
    \BEMAIL{elia.bisi@unifi.it}}\orcid{0000-0001-7628-8766}
  \and 
  Fabio~Deelan~Cunden\footnote{Universit\`a degli Studi di Bari and INFN, Sezione di Bari, Italy. \BEMAIL{fabio.cunden@uniba.it}}\orcid{0000-0002-9208-6467}}


\AMSSUBJ{60B20; 
60F15; 
62E15; 
44A60; 
05C05; 
05A15
} 
\AMSSUBJSECONDARY{33C20; 
46L54
} 

\VOLUME{0}
\YEAR{2025}
\PAPERNUM{0}
\DOI{10.1214/YY-TN}

\ABSTRACT{We consider random matrices whose shape is the dilation $N\lambda$ of a self-conjugate Young diagram $\lambda$.
In the large-$N$ limit, the empirical distribution of the squared singular values converges almost surely to a probability distribution $F^{\lambda}$.
The moments of $F^{\lambda}$ enumerate two combinatorial objects: $\lambda$-plane trees and $\lambda$-Dyck paths, which we introduce and show to be in bijection.
We also prove that the distribution $F^{\lambda}$ is algebraic, in the sense of Rao and Edelman.
In the case of fat hook shapes we provide explicit formulae for $F^{\lambda}$ and we express it as a free convolution of two measures involving a Marchenko--Pastur and a Bernoulli distribution.}

\definecolor{gray1}{gray}{0.9}
\definecolor{gray2}{gray}{0.75}
\definecolor{gray3}{gray}{0.6}

\newcommand{\R}{{\mathbb R}}
\newcommand{\C}{{\mathbb C}}
\newcommand{\N}{{\mathbb N}}
\newcommand{\Z}{{\mathbb Z}}
\newcommand\EE{{\mathbb E}}
\newcommand\Var{\operatorname{Var}}
\newcommand\Par{{\mathcal P}}
\newcommand\de{ \mathrm{d}}

\newcommand\Tr{\mathrm{Tr}}
\newcommand\tr{\mathrm{tr}}

\newcommand{\be}{\begin{equation}}
\newcommand{\ee}{\end{equation}}

\begin{document}



\section{Introduction}

In this article we study certain classes of random matrices associated with Young diagrams.
Let $\lambda$ be a Young diagram (i.e.\ integer partition), viewed as a set of boxes $(i,j)\in\N^2$.
A \emph{$\lambda$-shaped random matrix} ($\lambda$-RM for short), as defined in~\cite{Cunden23}, is a matrix $X$ such that the $(i,j)$-entries indexed by $(i,j)\in\lambda$ are i.i.d.\ complex random variables, while the other entries are zero.
We aim to study the empirical spectral distribution of random matrices of the form $XX^*$, where $X$ is a $\lambda$-RM.

When $\lambda$ is a rectangular Young diagram, $X$ simply has i.i.d.\ complex entries and $XX^*$ is a classical \emph{random covariance matrix}.
These matrix ensembles and their properties have been extensively studied in mathematics, statistics, and physics~\cite{Muirhead82,Forrester10,Livan18}.

In the special case of Gaussian entries, $\lambda$-RMs appeared in relation to representation theory, biorthogonal ensembles, last passage percolation, and free probability~\cite{Dykema04,Feray11,Adler13,Forrester17,Cheliotis18,Nakashima22,Forrester23}.

Recently, Cunden, Ligab\`o and Monni~\cite{Cunden23} considered $(N\lambda)$-RMs, where $N\lambda$ is the \emph{dilation} of a fixed Young diagram $\lambda$ (i.e., the Young diagram obtained from $\lambda$ by replacing every box with an $N\times N$ box grid).
In the case where $\lambda$ is a \emph{staircase partition}, they characterised the limiting spectral measure as $N\to\infty$, which is a generalisation of both the Marchenko--Pastur and Dykema--Haagerup measures.
Their analysis was based on a moment formula involving certain labelled trees (\emph{$r$-plane trees}), which were enumerated by Gu, Prodinger and Wagner~\cite{Gu10} in terms of a generalisation of Catalan numbers.
See Section~\ref{sec:preliminaries} for more details.

As observed in~\cite{Cunden23}, ``For $\lambda$-shaped random matrices one expects a relation between the limit shape $\lambda$ and the limiting spectral distribution.''
In the present work, we continue the research programme set up in~\cite{Cunden23}, investigating ``new examples, for which the limit [...]\ can be computed explicitly''.
In particular, we consider $(N\lambda)$-RMs, where $\lambda$ is a generic self-conjugate partition.
In Section~\ref{sec:mainResult}, we show that the empirical spectral distribution converges almost surely as $N\to\infty$ to a limit, whose moments are given in terms of enumerations of certain $\lambda$-dependent classes of labelled trees, which we call \emph{$\lambda$-plane trees} (see Theorem~\ref{thm:main}).
We also provide summation formulae for enumerating $\lambda$-plane trees (Theorem~\ref{thm:summationFormula}), which, in some cases, essentially reduce to terminating hypergeometric series (Corollary~\ref{coro:fat_hook}); a key ingredient for these results is a recent refined enumeration formula for $r$-plane trees of Okoth and Wagner~\cite{Okoth24}.
In Section~\ref{sec:lambdaDyckPaths}, we present a combinatorial bijection between $\lambda$-plane trees and certain $3$-dimensional lattice paths associated with the Young diagram $\lambda$, which we call \emph{$\lambda$-Dyck paths}.
In Section~\ref{sec:generatingFunctions}, we adopt a more analytic point of view, considering the generating functions of $\lambda$-plane trees and studying the systems of equations that they satisfy.
We show that $\lambda$-RMs are algebraic, in the sense that the Cauchy transform of their limiting spectral measure satisfies a polynomial equation (Theorem~\ref{thm:Stieltjes}).
Finally, we obtain an explicit solution for \emph{fat hooks}, i.e.\ Young diagrams made of two rectangular blocks: in this case, the limiting spectral measure, through its $R$-transform, can be  expressed as a free convolution of two measures involving a Marchenko--Pastur and a Bernoulli distribution (Theorem~\ref{thm:fatHook}).

\section{Preliminaries and background}
\label{sec:preliminaries}

In this section we start by recalling some standard notation related to integer partitions and introducing terminology that will be useful in the rest of the article.
Thereafter, we define the model of $\lambda$-RMs, recalling the result of~\cite{Cunden23} on the limiting spectral measure of random matrices associated with dilations of staircase partitions; this will be also our starting point in the next section.

Throughout, we will use the notation $\llbracket n \rrbracket :=\{1,\dots,n\}$.

An \textbf{(integer) partition} is a sequence $\lambda=\left(\lambda_1,\lambda_2,\ldots\right)$ of nonnegative integers, called \textbf{parts}, such that $\lambda_1\geq \lambda_2\geq \cdots$ and $\lambda_i=0$ for some $i\geq 1$.
Let $\Par$ be the set of all partitions.
For a given $\lambda\in\Par$, we let 
\[
|\lambda|:=\sum_{i\geq1}\lambda_i
\]
be its \textbf{size} and
\[
\ell(\lambda):=\min\{i\geq 0\colon \lambda_{i+1}=0\}
\]
be its \textbf{length}.

\begin{figure}
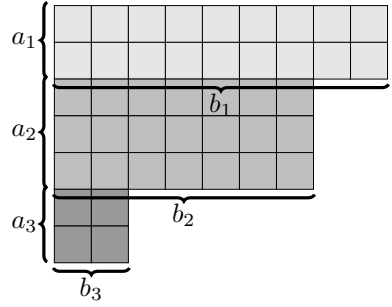

\centering
\[
\begin{ytableau}[*(gray2)]
        *(gray1)\tikznode{a1}{~}&*(gray1)&*(gray1)&*(gray1)&*(gray1)&*(gray1)&*(gray1)&*(gray1)&*(gray1)\\
        *(gray1)\tikznode{b11}{~}&*(gray1)&*(gray1)&*(gray1)&*(gray1)&*(gray1)&*(gray1)&*(gray1)&*(gray1)\tikznode{b12}{~}\\
        \tikznode{a2}{~}&&&&&&\\
        &&&&&&\\
        \tikznode{b21}{~}&&&&&&\tikznode{b22}{~}\\
        *(gray3)\tikznode{a3}{~}&*(gray3)\\
        *(gray3)\tikznode{b31}{~}&*(gray3)\tikznode{b32}{~}
\end{ytableau}
\]
\tikz[overlay,remember picture]{%
\draw[decorate,decoration={brace},very thick] ([yshift=-2mm,xshift=2mm]b12.south east) -- 
([yshift=-2mm,xshift=-2mm]b11.south west) node[midway,below]{$b_1$};
\draw[decorate,decoration={brace},very thick] ([yshift=-2mm,xshift=2mm]b22.south east) -- 
([yshift=-2mm,xshift=-2mm]b21.south west) node[midway,below]{$b_2$};
\draw[decorate,decoration={brace},very thick] ([yshift=-2mm,xshift=2mm]b32.south east) -- 
([yshift=-2mm,xshift=-2mm]b31.south west) node[midway,below]{$b_3$};
\draw[decorate,decoration={brace},very thick] ([yshift=-1.2mm,xshift=-2.8mm]b11.south west) -- 
([yshift=3.5mm,xshift=-2.8mm]a1.north west) node[midway,left]{$a_1$};
\draw[decorate,decoration={brace},very thick] ([yshift=-1.2mm,xshift=-2.8mm]b21.south west) -- 
([yshift=3.5mm,xshift=-2.8mm]a2.north west) node[midway,left]{$a_2$};
\draw[decorate,decoration={brace},very thick] ([yshift=-1.2mm,xshift=-2.8mm]b31.south west) -- 
([yshift=3.5mm,xshift=-2.8mm]a3.north west) node[midway,left]{$a_3$};
}
\caption{The Young diagram corresponding to the partition $\lambda=(9,9,7,7,7,2,2)$ of length $\ell(\lambda)=7$ and size $|\lambda|=43$.
It has $3$ rectangular blocks, depicted in different shades of gray, with respective heights $\bm{a}=(a_1,a_2,a_3)=(2,3,2)$ and bases $\bm{b}=(b_1,b_2,b_3)=(9,7,2)$.
In multirectangular coordinates, $\lambda$ can thus be written as $\bm{a}\times \bm{b}=(9^2,7^3,2^2)$.
The block map associated with $\lambda$ is $f_{\lambda}\colon \llbracket 7\rrbracket \to \llbracket 3\rrbracket$ with $f_{\lambda}(1)=f_{\lambda}(2)=1$, $f_{\lambda}(3)=f_{\lambda}(4)=f_{\lambda}(5)=2$, and $f_{\lambda}(6)=f_{\lambda}(7)=3$.
}
\label{fig:partition}
\end{figure}

We identify a partition $\lambda\in\Par$ with its \textbf{Young diagram}, i.e.\ the set
\be
\left\{(i,j)\in\Z^2\colon i\in \llbracket \ell(\lambda)\rrbracket, j\in \llbracket \lambda_i\rrbracket\right\},
\ee
which is usually depicted as a set of left-aligned boxes, where the $i^{\text{th}}$ row contains $\lambda_i$ boxes (see Fig.~\ref{fig:partition}).
Under the involution $(i,j)\mapsto(j,i)$ of $\Z^2$, the Young diagram of $\lambda$ is mapped into the Young diagram of another partition $\lambda'\in\Par$, called the \textbf{conjugate partition}.
If $\lambda=\lambda'$, we say that $\lambda$ is \textbf{self-conjugate} (see Fig.~\ref{fig:lambda-plane-tree}).

It is sometimes  convenient to use \textbf{multirectangular coordinates} (see~\cite{Stanley06}).
For two lists of positive integers $\bm{a}=(a_1,\ldots,a_r)$ and $\bm{b}=(b_1,\ldots,b_r)$, with $b_1> b_2>\cdots> b_r$, we define the partition 
\be
\bm{a}\times\bm{b}:=(\underbrace{b_1,\ldots,b_1}_{\text{$a_1$ times}},\underbrace{b_2,\ldots,b_2}_{\text{$a_2$ times}},\ldots,\underbrace{b_r,\ldots,b_r}_{\text{$a_r$ times}})=:(b_1^{a_1},b_2^{a_2},\ldots,b_r^{a_r})\in\Par,
\ee
which is of length $\sum_{i=1}^r a_i$ and size $\sum_{i=1}^r a_i b_i$.
The corresponding Young diagram is a union of $r$ rectangular blocks of sizes $a_i\times b_i$, where $a_i$ is the height and $b_i$ is the base, with strictly decreasing bases (see Fig.~\ref{fig:partition}). 
Any $\lambda\in\Par$ can be represented uniquely in multirectangular coordinates as follows: set $b_1:=\lambda_1$ and, recursively for all $i\geq 2$, $b_i:=\lambda_{\min\{j\geq 2\colon \lambda_j<b_{i-1}\}}$ (the sequence $b_1,b_2,\dots$ stops at $b_r$ if $\lambda$ has no positive parts smaller than $b_r$); finally, set $a_i$ to be the number of parts of $\lambda$ equal to $b_i$, for all $i\in\llbracket r \rrbracket$; then, we have $\lambda=\bm{a}\times\bm{b}$.

We call \textbf{block map} any non-decreasing, surjective function $f\colon \llbracket \ell \rrbracket\to\llbracket r \rrbracket$, for some positive integers $\ell\geq r$.
The \textbf{block map associated with $\lambda=\bm{a}\times \bm{b}$} is the function
\be
f_\lambda\colon \llbracket \ell(\lambda)\rrbracket\to\llbracket r \rrbracket, \qquad
f_{\lambda}(i)=n \quad \text{if} \quad \lambda_i=b_n.
\ee
Notice that, since $b_1> b_2>\cdots> b_r$, the condition $\lambda_i=b_j$ is satisfied for one and only one index $j$.
In words, $f_\lambda(i)=j$ if and only if the $i^{\text{th}}$ row of the Young diagram $\lambda$ belongs to the $j^{\text{th}}$ rectangular block in the multirectangular coordinate representation.
Notice that any block map $f\colon \llbracket \ell \rrbracket\to\llbracket r \rrbracket$ equals $f_\lambda$ for infinitely many choices of $\lambda\in\Par$ of length $\ell$ and with $r$ blocks: the reason is that a block map determines the heights $a_i$'s of the blocks, but not the their bases $b_i$'s.

If $\lambda=\bm{a}\times \bm{b}$, then $\lambda=\lambda'$ if and only if $b_i=\sum_{i=1}^{r-i+1}a_i$.
Thus, any self-conjugate partition can be written as
     \be
     \label{eq:selfConjMultirectangular}
     \lambda=((a_1+a_2+\cdots+a_r)^{a_1},(a_1+a_2+\cdots+a_{r-1})^{a_2},\ldots,a_1^{a_r})
     \ee
for some list $\bm{a}=(a_1,\ldots,a_r)$ of positive integers.
We can then also write
\be\label{eq:selfConjMultirectangular2}
\lambda_i =a_1+\cdots +a_{r-f_{\lambda}(i)+1}, \qquad
i\in \llbracket \ell(\lambda)\rrbracket.
\ee
If one restricts to self-conjugate partitions, there is a one-to-one correspondence between self-conjugate partitions and block maps, since, when $\lambda=\lambda'$, the $b_i$'s are determined by the $a_i$'s.
Under this correspondence, for instance, the staircase partition $\lambda=(r,r-1,\ldots,1)$ is associated with the identity map on $\llbracket r \rrbracket$.

Furthermore, we have the following simple characterisation:
\begin{lemma}
\label{lem:blockMap}
Let $\lambda=\lambda'\in \Par$ have $r$ blocks.
For all $i,j\in \llbracket \ell(\lambda)\rrbracket$, 
\be\label{eq:lambda_equiv}
(i,j)\in\lambda \qquad\Longleftrightarrow\qquad
f_{\lambda}\left(i\right)+f_{\lambda}\left(j\right)\leq r+1.
\ee
\end{lemma}
\begin{proof}
Let $\lambda$ as in~\eqref{eq:selfConjMultirectangular}.
If $(i,j)\in\lambda$, then, by definition, $j\leq \lambda_i$.
Using the monotonicity of $f_{\lambda}$ and equality~\eqref{eq:selfConjMultirectangular2}, we have
\be
\label{eq:blockMap}
f_{\lambda}(j)
\leq f_{\lambda}(\lambda_i)
= f_{\lambda}(a_1+\dots+a_{r-f_{\lambda}(i)+1})
= r - f_{\lambda}(i)+1 .
\ee

Conversely, if $f_{\lambda}\left(i\right)+f_{\lambda}\left(j\right)\leq r+1$, then 
\[
j\leq \max\{h\in \llbracket \ell(\lambda) \rrbracket \colon f_{\lambda}(h)\leq r-f_{\lambda}(i)+1\}
= a_1+\dots + a_{r-f_{\lambda}(i)+1}
= \lambda_i,
\]
where the latter equality follows again from~\eqref{eq:selfConjMultirectangular2}.
This proves that $(i,j)\in\lambda$, as claimed.
\end{proof}

There is a natural operation of multiplication of partitions by positive integers.
If $\lambda=(\lambda_1,\lambda_2,\ldots)$ and $N\in\N$, then the \textbf{dilation} of $\lambda$ by $N$ is the partition
\[
N\lambda:=(\underbrace{N\lambda_1,\ldots,N\lambda_1}_{\text{$N$ times}},\underbrace{N\lambda_2,\ldots,N\lambda_2}_{\text{$N$ times}},\ldots).
\]
We can also think of $N\lambda$ as the Young diagram obtained by replacing every box of $\lambda$ with an $N\times N$ grid of boxes.
We have $|N\lambda|=N^2|\lambda|$ and $\ell(N\lambda)=N\ell(\lambda)$. 

We say that $\lambda\in\Par$ is \textbf{minimal} if, for any $\mu\in\Par$,
\[
\lambda=N\mu \quad\iff\quad \text{$N=1$ and $\mu=\lambda$}.
\]
It is easy to see that $\lambda\in\Par$ is minimal if and only if $\gcd(\lambda_1,\lambda_2,\ldots,\lambda'_1,\lambda'_2,\dots)=1$; in multirectangular coordinates, $\lambda=\bm{a}\times\bm{b}$ is minimal if and only if $\gcd(a_1,\dots,a_r,b_1,\dots,b_r)=1$.
By~\eqref{eq:selfConjMultirectangular}, a self-conjugate partition $\lambda=\lambda'$ is minimal if and only if $\gcd(a_1,\ldots,a_r)=1$.

To set up the model of $\lambda$-shaped random matrices, let $\{X_{i,j}\colon i,j\in\N\}$ be a collection of i.i.d.\ complex random variables with $\EE X_{i,j}=0$ and $\EE |X_{i,j}|^2=1$.
Let $\lambda\in\Par$ be a fixed partition of length $\ell(\lambda)=:\ell$.
For $N\in\N$, let $X_N$ be the $(N\ell)\times (N\ell)$ random matrix whose $(i,j)$-entry is $X_{i,j}$ if $(i,j)\in N\lambda$, and $0$ otherwise.
Then, $X_N$ is an $(N\lambda)$-shaped random matrix, according to the terminology of the introduction.
Consider the $(N\ell)\times (N\ell)$ nonnegative definite complex Hermitian matrix
\be
\label{eq:defW}
W_N:=\frac{1}{N {\ell}} X_N X_N^*.
\ee
and denote by $x_1^{(N)}\leq x_2^{(N)}\leq\cdots\leq x_{N\ell}^{(N)}$ its eigenvalues, i.e.\ the squared singular values of $X_N/\sqrt{N{\ell}}$.
Let
\be
\label{eq:defF}
F_N(x):=\frac{1}{N\ell}\#\left\{j\in \llbracket N\ell\rrbracket \colon x_j^{(N)}\leq x \right\}, \qquad x\in\R,
\ee
be the empirical distribution (function) of the eigenvalues of $W_N$.
It is natural to ask whether the sequence of \textbf{spectral distributions} $(F_N)_{N\geq1}$ converges to some distribution.
The limit, when it exists, will be referred to as the \textbf{limiting spectral distribution} of $W_N$ and denoted $F^{\lambda}$.

It is clear that, if $\lambda=N\mu$ for some minimal partition $\mu$, then the sequence of spectral distributions associated with $\lambda$ is a subsequence of the corresponding sequence for $\mu$.
If the limit $F^{\mu}$ exists, then $F^{\lambda}$ also does {and $F^{\lambda}=F^{\mu}$}.
For this reason, we may suppose from now on that $\lambda$ is a minimal partition.

Note that, for any fixed $\lambda\in\Par$, the sequence of nested Young diagrams $(\lambda^{(N)})_{N\geq1}$ defined by $\lambda^{(N)}:=N\lambda$, is \textbf{balanced}~\cite{Dolega10}, in the sense that
{both the number of rows and the number of columns of $\lambda^{(N)}$ are $O\big(\sqrt{|\lambda^{(N)}|}\big)$} as $N\to\infty$.
Such a condition is known~\cite{Cunden23} to be a necessary condition for the existence of the limit $F^{\lambda}$.

A limit theorem in the case where $\lambda$ is a staircase partition was obtained in~\cite{Cunden23}.
  \begin{theorem}[\cite{Cunden23}]
  \label{thm:staircase}
  Let $\lambda=(r,r-1,\dots,1)\in\Par$.
  Then, almost surely, $F_N$ converges to the deterministic distribution $F^{(r,r-1,\dots,1)}$ on $[0,\infty)$ whose moments are
  \be
  \label{eq:staircase}
  \int_0^{\infty} x^k \,\de F^{(r,r-1,\dots,1)}(x)=\frac{1}{{r^k(k+1)}} \binom{(r+1)k}{k}, \qquad k\geq 0.
  \ee
  \end{theorem}
  {Notice that the normalisation in~\eqref{eq:defW} differs from the one considered in~\cite{Cunden23}, resulting in a slightly different expression for the moments~\eqref{eq:staircase}.}
A key ingredient in the proof of the latter theorem was the following combinatorial object introduced in~\cite{Gu10}.
\begin{definition}
\label{def:rPlaneTrees}
Let $r\geq1$.
An \textbf{$r$-plane tree} is a pair $(T,c)$, where $T$ is a plane tree with vertex set $V$ and edge set $E$, and $c\colon V \to\llbracket r \rrbracket$ is a labelling such that $c(u)+c(v)\leq r+1$ whenever $\{u,v\}\in E$.
\end{definition}
It turned out that the limiting moments on the left-hand side of~\eqref{eq:staircase} are related to enumerations of $r$-plane trees, hence Theorem~\ref{thm:staircase} was proved using the following formula.
    \begin{theorem}[\cite{Gu10}]
    \label{thm:GPW}
The number of $r$-plane trees on $k+1$ vertices is
  \be
  \label{eq:CatalanGen}
\frac{r}{k+1}\binom{(r+1)k}{k}, \qquad k\geq 0.
  \ee
  \end{theorem}
  
For $r=1$, so that $\lambda=(1)$ is the partition with one box only, \eqref{eq:CatalanGen} is the sequence of Catalan numbers and the limiting distribution $F^{(1)}$ in Theorem~\ref{thm:staircase} is the classical Marchenko--Pastur law.

{Notice that the integer sequence~\eqref{eq:CatalanGen} resembles the Fuss--Catalan numbers, which arise as limiting moments of products of random matrices; see e.g.~\cite{AlexeevGoetzeEtAl10, Mlotkowski10, LiuSongWang11, PensonZyczkowski11, BanicaBelinschiEtAl11}.}
  
Gu, Prodinger and Wagner~\cite{Gu10} concluded their article with the following remark (notation adapted to the present work):
``\emph{one can certainly modify the definition of $r$-plane trees by imposing other restrictions on pairs of labels along an edge. It is conceivable that appropriate conditions will lead to interesting counting problems as well.}''.
As we will see in Section~\ref{sec:mainResult}, it is precisely a more general condition on pairs of edges in Definition~\ref{def:rPlaneTrees} that lead us to characterise the limiting spectral measure of $\lambda$-RMs, with $\lambda$ any self-conjugate partition.

\section{$\lambda$-shaped random matrices and $\lambda$-plane trees}
\label{sec:mainResult}

In this section we state and prove a limit theorem for the spectral distribution of $(N\lambda)$-shaped random matrices, where $\lambda$ is a self-conjugate Young diagram.
The moments of the limiting spectral distribution are given in terms of enumerations of a combinatorial object, which generalises $r$-plane trees, and for which we provide a summation formula.

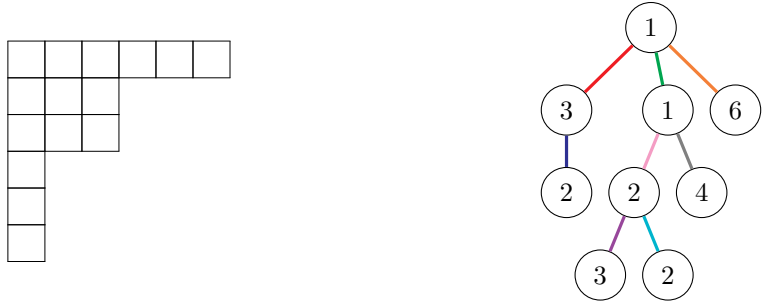
\begin{figure}
\begin{subfigure}[c]{0.5\textwidth}
\centering
\ydiagram{6,3,3,1,1,1}
\end{subfigure}%
\begin{subfigure}[c]{0.5\textwidth}
\centering   
 \begin{forest}
 for tree={edge={very thick}}
[$1$
[$3$,edge=Red
[$2$,edge=Blue
]]
[$1$,edge=Green
[$2$,edge=Lavender
[$3$,edge=Purple
]
[$2$,edge=Turquoise
]]
[$4$,edge=gray
]]
[$6$,edge=Orange
]]
\end{forest}
\end{subfigure}%
\caption{On the left-hand side, the self-conjugate partition $\lambda=\lambda'=(6^1,3^2,1^3)$ of length $\ell(\lambda)=6$.
On the right-hand side, a $\lambda$-plane tree: every vertex has a label in $\llbracket6\rrbracket$ and every edge $\{u,v\}$ satisfies $(c(u),c(v))\in\lambda$.
The edges of the tree are depicted in different colours for later purposes.}
\label{fig:lambda-plane-tree}
\end{figure}

\begin{definition}
\label{def:lambda-plane}  Let $\lambda\in\Par$ with $\lambda=\lambda'$ and $\ell(\lambda)=:\ell$.
A \textbf{$\lambda$-plane tree} is a pair $(T,c)$, where $T$ is a plane tree with vertex set $V$ and edge set $E$, and $c\colon V \to\llbracket \ell\rrbracket$ is a labelling such that $(c(u),c(v))\in\lambda$ whenever $\{u,v\}\in E$.
See Figure~\ref{fig:lambda-plane-tree}.
\end{definition}

By Lemma~\ref{lem:blockMap}, we can rewrite the condition on the labelling $c$ of Definition~\ref{def:lambda-plane} in terms of the block map $f_{\lambda}\colon \llbracket \ell\rrbracket \to \llbracket r\rrbracket$ associated with $\lambda$:
\be\label{eq:lambdaPlaneTrees_equiv}
(c(u),c(v))\in\lambda \qquad\Longleftrightarrow\qquad
f_{\lambda}\left(c(u)\right)+f_{\lambda}\left(c(v)\right)\leq r+1.
\ee

\begin{example} If $\lambda=\lambda'=\bm{a}\times \bm{b}$ with $\bm{a}=(a_1,\ldots,a_r)=(1,\ldots,1)$, then $\bm{b}=(r,r-1,\dots,1)$, so that $\lambda=(r,r-1,\ldots,1)$ is the staircase partition of length $\ell=r$ and its block map is the identity on $\llbracket r\rrbracket$.
It then follows from~\eqref{eq:lambdaPlaneTrees_equiv} that the corresponding $\lambda$-plane trees are precisely the $r$-plane trees of Definition~\ref{def:rPlaneTrees}.
\end{example}

\begin{example}
For $r=2$ and $\bm{a}=(a_1,a_2)$, we call $\lambda=((a_1+a_2)^{a_1},a_1^{a_2})$ a \textbf{fat hook}.
The condition on the labelling of the trees is that $\min(c(u),c(v))\leq a_1$ whenever $u$ and $v$ are adjacent vertices.
Special cases include the following.
Fix $\ell\geq2$, and set $a_1=1$ and $a_2=\ell-1$, then $\lambda=(\ell^1,1^{\ell-1})$ is a \textbf{hook}; in this case, the block map is $f_{\lambda}(i)=\min(i,2)$, and the condition on the labelling of the trees is that $\min(c(u),c(v))= 1$ whenever $u$ and $v$ are adjacent vertices.
If, instead, we set $a_1=\ell-1$ and $a_2=1$, then $\lambda=(\ell^{\ell-1},(\ell-1)^1)$ is a \textbf{notched square}; in this case, the block map is $f_{\lambda}(i)=\max(i-\ell+2,1)$, and the condition on the labelling of the trees is that $\min(c(u),c(v))\leq \ell-1$ whenever $u$ and $v$ are adjacent vertices.
 \end{example}

We are now ready to present our limit theorem, which is based on the moment method: we characterise the limiting spectral measure in terms of its moments, which in turn are expressed in terms of the integer sequence
\be
\label{eq:defClk}
  C_{k}^{\lambda}:=\#\{\text{$\lambda$-plane trees on $k+1$ vertices}\}, \qquad k \geq 0.
\ee
  \begin{theorem}
  \label{thm:main}   Let $\lambda=\lambda'\in\Par$ with $\ell(\lambda)=:\ell$.
  Then, as $N\to\infty$,  the sequence of random distributions $F_N$ defined in~\eqref{eq:defF} converges almost surely to the deterministic distribution $F^{\lambda}$ on $[0,\infty)$ whose moments are 
  \be\label{eq:limitMoments}
  \int_0^{\infty} x^k \,\de F^{\lambda}(x)=\frac{C^{\lambda}_k}{\ell^{k+1}}. 
  \ee
  \end{theorem}

  When $r=\ell$ and $\lambda=(r,r-1,\dots,1)$ is a staircase partition, we recover Theorem~\ref{thm:staircase}; in this case, a closed formula for $C_{k}^{\lambda}$ is given in Theorem~\ref{thm:GPW}.
  In the general case, we obtain a summation formula for this combinatorial sequence.
          \begin{theorem}
          \label{thm:summationFormula} Let $\lambda=\lambda'\in\Par$, written in multirectangular coordinates as
        \[
        \lambda=((a_1+a_2+\cdots+a_r)^{a_1},(a_1+a_2+\cdots+a_{r-1})^{a_2},\ldots,(a_1)^{a_r}).
        \]
         Then, the number of $\lambda$-plane trees on $k+1$ vertices is
  \be
  \label{eq:summationFormula}
  C_{k}^{\lambda}=\sum_{\substack{\ell_1,\ldots,\ell_r\geq 0,\\ \ell_1+\cdots+\ell_r=k+1}} t(\ell_1,\ldots,\ell_r) \, a_1^{\ell_1}\cdots a_r^{\ell_r},
\ee
where, using the shorthand $\ell_{\leq i}:=\ell_1+\cdots+\ell_i$ and $\ell_{\geq i}:=\ell_{i}+\cdots+\ell_r$ for any $i\in\llbracket r\rrbracket$, we set
\be\label{eq:binomProd}
t(\ell_1,\ldots,\ell_r)
:= \frac{1}{k}\prod_{j=1}^{r}\binom{\ell_{\geq j}+\ell_{\leq r-j+1}-1-{{\bm{1}}_{\llbracket \lceil{r/2}\rceil \rrbracket}(j)}}{\ell_j}.
\ee
  \end{theorem}
The proof of Theorem~\ref{thm:main} and~\ref{thm:summationFormula} is postponed to the end of this section.
We first look at the case of fat hooks, where the summation in~\eqref{eq:summationFormula} essentially reduces to a terminating hypergeometric series.
      \begin{corollary}\label{coro:fat_hook} Let $\lambda=((a_1+a_2)^{a_1},a_1^{a_2})$ be a fat hook.
      Then,
  \[
  C_{k}^{\lambda}
  =C_k \, a_1^{k+1} \, {}_{2}F_{1} \left( \begin{matrix}
 -k-1,-k \\ 
k \end{matrix} ;\frac{a_2}{a_1} \right),
\]
where $C_k:=\frac{1}{k+1}\binom{2k}{k}$ is the $k^{\text{th}}$ Catalan number and
\be
{}_{2}F_{1} \left( \begin{matrix}
 \alpha,\beta \\ 
\gamma \end{matrix} ;x \right)
:= \sum_{j=0}^{\infty} \frac{(\alpha)_j (\beta)_j}{(\gamma)_j } \frac{x^j}{j!}
\ee
is the hypergeometric function (here, $(\alpha)_j:= \alpha(\alpha+1)\cdots(\alpha+j-1)$ denotes the rising Pochhammer symbol).
  \end{corollary}
  
  \begin{proof}
From~\eqref{eq:binomProd}, we have   
\[
 t(\ell_1,\ell_2)= \frac{1}{k}\binom{2k}{\ell_1}\binom{k}{\ell_2}  \qquad \text{if} \;\; \ell_1,\ell_2\geq 0, \;\; \ell_1+ \ell_2= k+1.
\]
Inserting this into~\eqref{eq:summationFormula}, we obtain
  \[
  \begin{aligned}
    C_{k}^{\lambda}&=\sum_{\substack{\ell_1,\ell_2\geq 0\\ \ell_1+\ell_2=k+1}} \frac{1}{k}\binom{2k}{\ell_1}\binom{k}{\ell_2} a_1^{\ell_1}a_2^{\ell_2}
    =
    \sum_{j=0}^{k+1} \frac{1}{k}\binom{2k}{k+1-j}\binom{k}{j}a_1^{k+1-j}a_2^{j} \\
&=   \frac{1}{k+1}\binom{2k}{k}    a_1^{k+1} \sum_{j=0}^{k+1} \frac{1}{j!} \frac{(k+1)!}{(k+1-j)!} \frac{(k)!}{(k-j)!} \frac{(k-1)!}{(k-1+j)!}\left(\frac{a_2}{a_1}\right)^{j} \\
&=C_k \, a_1^{k+1} \, {}_{2}F_{1} \left( \begin{matrix}
 -k-1,-k \\ 
k \end{matrix} ;\frac{a_2}{a_1} \right),
 \end{aligned}
 \]
as claimed.
  \end{proof}

\begin{example}
Two special cases of Corollary~\ref{coro:fat_hook} are the following:
      \begin{enumerate}
    \item
If $\lambda=(\ell,1^{\ell-1})$ is a hook, then
  \[
  C_{k}^{(\ell,1^{\ell-1})}=C_k \, {}_{2}F_{1} \left( \begin{matrix}
 -k-1,-k \\ 
k \end{matrix} ;\ell-1 \right).
\]
\item
If $\lambda=(\ell^{\ell-1},\ell-1)$ is a notched square, then
  \[
  C_{k}^{(\ell^{\ell-1},\ell-1)} = C_k \, (\ell-1)^{k+1} \, {}_{2}F_{1} \left( \begin{matrix}
 -k-1,-k \\ 
k \end{matrix} ;\frac{1}{\ell-1} \right).
  \]
  \end{enumerate}
\end{example}

  \begin{remark}
Let $X$ be a matrix with i.i.d.\ standard complex \emph{Gaussian} entries $X_{i,j}$ for $(i,j)\in\lambda$ and $X_{i,j}=0$ for $(i,j)\notin\lambda$.
Thus, adapting the terminology of the introduction, $X$ is a Gaussian $\lambda$-RM.
In such a Gaussian case, a combinatorial description of the moments of the matrix $XX^*$ was obtained in~\cite[Theorem 7]{Feray11} using Wick's formula:
 \be
 \label{eq:FeraySnaidy}
 \EE \Tr \left(XX^*\right)^k=\sum_{\substack{\sigma_1,\sigma_2\in S_k\colon \\ \sigma_1\circ \sigma_2=\gamma_k}}\mathrm{N}^{\lambda}(\sigma_1,\sigma_2),
 \ee
 where  $\gamma_k=(1,\ldots,k)$ is the full cycle in the symmetric group $S_k$, and $\mathrm{N}^{\lambda}(\sigma_1,\sigma_2)$ denotes the number of colourings of the cycles of $\sigma_1$ and $\sigma_2$ that are compatible with $\lambda$ (see~\cite{Feray11} for the precise definition).
Theorems~\ref{thm:main}-\ref{thm:summationFormula} give the leading term of the large-$N$ asymptotics of formula~\eqref{eq:FeraySnaidy} when the shape is $N\lambda$ with $\lambda$ self-conjugate.
  \end{remark}
  
  We now prove Theorems~\ref{thm:main} and~\ref{thm:summationFormula}. 
  \begin{proof}[Proof of Theorem~\ref{thm:main}]
The proof generalises the argument presented in~\cite{Cunden23} and is based on the moment method for random matrices.
We first check that the distribution $F^{\lambda}$ is uniquely determined by the moment sequence $m_k:= C_{k}^{\lambda}/{\ell^{k+1}}$.
To do this, we check the Riesz condition~\cite[Lemma~B.2]{Bai10}:
  \[
  \liminf_k\frac{1}{k}m_{2k}^{\frac{1}{2k}}<\infty.
  \]
  Notice that the number of plane trees on $k+1$ vertices where each vertex is labelled with \emph{any} integer in $\llbracket \ell\rrbracket$ equals $\ell^{k+1}C_{k}$, where $C_{k}=\frac{1}{k+1}\binom{2k}{k}$ are the Catalan numbers.
  Therefore, 
  \be\label{eq:lambdaPlaneNoBound}
  C^{\lambda}_{k}\leq \ell^{k+1}C_{k}
  \ee
  and
   \[
  \frac{1}{k}m_{2k}^{\frac{1}{2k}}
  =\frac{1}{k}\left[\frac{1}{\ell^{2k+1}}C^{\lambda}_{2k}\right]^{\frac{1}{2k}}
  \leq \frac{1}{k}{C_{2k}^{\frac{1}{2k}}}
  =\frac{1}{k}\left[\frac{1}{2k+1}\binom{4k}{2k}\right]^{\frac{1}{2k}}
  \xrightarrow{k\to\infty} 0,
  \]
  since {$\binom{4k}{2k}\leq 2^{4k}$}.
  
Consider now the sequence
\[
m_{k,N}:=\int_0^{\infty} x^k \,\de F_N(x)  =\frac{1}{N\ell} \Tr W_N^k
\]
 of moments of $F_N$.
  By the Moment Convergence Theorem~\cite[Section B.1]{Bai10}, the claim will follow if we prove that, for all integers $k\geq0$, 
  \[
m_{k,N}
  \xrightarrow{N\to\infty} m_k\qquad \text{a.s.}.
  \]
  To prove the latter, using Borel-Cantelli lemma and Chebyshev's inequality, it is enough to show that, for all integers $k\geq0$,
  \begin{align}
  \label{eq:lem1}
  \lim_{N\to\infty}\EE m_{k,N}&=m_k \qquad\text{and}\\
    \label{eq:lem2}
   \sum_{N=1}^{\infty} \Var m_{k,N}&<\infty.
  \end{align}
Using definition~\eqref{eq:defW}, we compute
\[
\begin{split}
&m_{k,N}
=\frac{1}{N\ell}\Tr W_N^k
=\frac{1}{(N\ell)^{k+1}}\sum_{i,j\colon \llbracket k\rrbracket \to \llbracket N\ell\rrbracket}\prod_{m=1}^{k}(X_N)_{i(m),j(m)}\overline{ (X_N)_{i(m+1),j(m)}},
\end{split}
\]
using the convention $i({k+1})\equiv i(1)$, which will be always assumed to hold from now on.
Recalling that, by definition, $(X_N)_{p,q}=X_{p,q} {\bm{1}_{N\lambda}(p,q)}$ for all $(p,q)\in {\llbracket N\ell\rrbracket \times \llbracket N\ell\rrbracket}$, we have
\be
\label{eq:moment}
m_{k,N} = \frac{1}{(N\ell)^{k+1}}  \sum_{(i,j)\in \Lambda^{(N)}} \prod_{m=1}^k X_{i(m),j(m)} \overline{X_{i(m+1),j(m)}},
\ee
where
\[
\Lambda^{(N)}:= \left\{(i,j) \text{ such that } i,j\colon \llbracket k\rrbracket \to \llbracket N\ell\rrbracket, \, (i(m),j(m))\in N\lambda, \, (i(m+1),j(m))\in N\lambda\right\}.
\]

We now re-express the above formula in terms of certain combinatorial objects.
For any graph $G=(V,E)$, we define a \textbf{walk} on $G$ to be a (finite or infinite) sequence $w=(w_0,w_1,\dots)$ in $V$ such that $\{w_i,w_{i+1}\}\in E$ for all $i\geq 0$.
We say that $w$ is \textbf{alternating} if $\{w_0,w_2,\dots\}\cap\{w_1,w_3,\dots\}=\emptyset$, and \textbf{spanning} if $\{w_0,w_1,\dots\}=V$.
We also say that a finite walk $ w=(w_0,\dots,w_{2k})$ of length $2k$ is \textbf{closed} if $w_0=w_{2k}$.
Thus, an alternating spanning closed walk (ASCW) is a walk on the graph that starts and ends at the same vertex and visits each vertex at least once, alternating visits to vertices from some subset of $V$ and visits to vertices from its complement.

Given $(i,j)\in \Lambda^{(N)}$, we construct a graph $G(i,j)$ with vertex set
\be
\label{eq:vertexSet}
V(i,j) := \{(1,i(m))\colon m\in \llbracket k\rrbracket\} \cup \{(2,j(m))\colon m\in \llbracket k\rrbracket\}
\ee
and edge set
\[
E(i,j) := \{\{(1,i(m)), (2,j(m))\}\colon m\in \llbracket k\rrbracket\} \cup \{\{(2,j(m)), (1,i(m+1))\}\colon m\in \llbracket k\rrbracket\},
\]
where both $V(i,j)$ and $E(i,j)$ might have less than $2k$ elements because of possible repetitions.
These definitions specify also a natural \textbf{labelling} of the vertices of $G(i,j)$ defined by
\be
\label{eq:labellingN}
\mathfrak{c}\colon V(i,j)\to \llbracket N\ell \rrbracket, \qquad
\mathfrak{c}((1,p))=\mathfrak{c}((2,p)):= p \qquad
\forall p\in\llbracket N\ell \rrbracket.
\ee
Moreover, the sequence
\[
 w(i,j) := (\underbrace{(1,i(1))}_{=w_0(i,j)}, \underbrace{(2,j(1))}_{=w_1(i,j)}, \underbrace{(1,i(2))}_{=w_2(i,j)}, \dots, \underbrace{(1,i(k))}_{=w_{2k-2}(i,j)}, \underbrace{(2,j(k))}_{=w_{2k-1}(i,j)}, \underbrace{(1,i(1))}_{=w_{2k}(i,j)})
\]
is an ASCW of length $2k$ on $G(i,j)$.
We may thus rewrite~\eqref{eq:moment} as
\be
\label{eq:moment2}
m_{k,N} = \frac{1}{(N\ell)^{k+1}} \sum_{(i,j)\in \Lambda^{(N)}} X_{ w(i,j)} ,
\ee
where $X_{ w(i,j)}$ is the product of random variables appearing in~\eqref{eq:moment} and determined by $ w(i,j)$:
\be
\label{eq:prodRVsGraph}
X_{ w(i,j)} := \prod_{m=1}^k X_{\mathfrak{c}(w_{2m-2}(i,j)), \mathfrak{c}(w_{2m-1}(i,j))} \overline{X_{\mathfrak{c}(w_{2m}(i,j)), \mathfrak{c}(w_{2m-1}(i,j))}} .
\ee

With this notation, we can proceed to show~\eqref{eq:lem1}-\eqref{eq:lem2}.
We may assume that there exists a positive constant $M$ such that $|X_{p,q}|\leq M$ almost surely.
{Such a boundedness assumption can be lifted by adapting the standard truncation argument explained in~\cite[p.~48]{Bai10}, if we note that 
\[
\frac{1}{\ell N}\operatorname{rank} \EE\widehat{X}_N\leq \frac{r}{\ell N}\xrightarrow{N\to\infty}0,
\]
where $\widehat{X}_N$ is the truncated matrix and $r$ is the number of blocks of $\lambda$.}

To prove~\eqref{eq:lem1}, start with
\be
\label{eq:momentMean}
\EE m_{k,N}
=\frac{1}{(N\ell)^{k+1}} \sum_{(i,j)\in \Lambda^{(N)}} \EE X_{ w(i,j)} .
\ee
If $ w(i,j)$ traverses a given edge only once, then $\EE X_{ w(i,j)}=0$, because the $X_{p,q}$'s are independent and centred.
Thus, all nonzero summands in~\eqref{eq:momentMean} correspond to $(i,j)$ such that the graph $G(i,j)$ has at most $k$ edges, and so at most $k+1$ vertices.
For $(i,j),(i',j')\in \Lambda^{(N)}$, we say that $ w(i,j)$ is isomorphic to $ w(i',j')$ if, up to renaming the vertices, they are the same ASCWs on the same graph.
For $n\leq k+1$, let $ W_n$ be the set of (isomorphism classes of) ASCWs of length $2k$ on some graph with \emph{exactly} $n$ vertices and \emph{at most} $k$ edges.
If $ w(i,j)$ is in the isomorphism class $ w\in  W_n$, we write $[ w(i,j)]= w$.
We thus have 
\[
\EE m_{k,N}
=\frac{1}{(N\ell)^{k+1}} \sum_{n=1}^{k+1} \sum_{ w\in  W_n}  \sum_{\substack{(i,j)\in \Lambda^{(N)}\colon \\ [ w(i,j)]= w}} \EE X_{ w(i,j)} .
\]

Let us now fix an ASCW $ w\in W_n$ on a graph $G$ with $n$ vertices.
The number of $(i,j)\in\Lambda^{(N)}$ such that $[ w(i,j)]= w$ is certainly bounded above by $(N\ell)^n$, and we have $\EE |X_{ w(i,j)}|\leq M^{2k}$.
It follows that the only nonzero contribution for large $N$ corresponds to $n=k+1$:
\[
\lim_{N\to\infty} \EE m_{k,N}
=\lim_{N\to\infty} \frac{1}{(N\ell)^{k+1}} \sum_{ w\in  W_{k+1}}  \sum_{\substack{(i,j)\in \Lambda^{(N)}\colon \\ [ w(i,j)]= w}} \EE X_{ w(i,j)} .
\]

Fix now $w\in W_{k+1}$.
The graph $G$ on which the walk $w$ travels has exactly $k+1$ vertices and at most $k$ edges, hence it is a tree with exactly $k$ edges.
Since $w$ is a closed spanning walk of length $2k$, it traverses each edge of $G$ exactly twice, in opposite directions.
Thus, $w$ may be identified with a (rooted) plane tree with $k+1$ vertices (see the beginning of the proof of Theorem~\ref{thm:lambdaDyck} for more details on how the sequence $ w$ encodes, in a standard way, the root and the `contour' or planar structure of a plane tree).
Moreover, since the $X_{p,q}$'s are independent with $\EE |X_{p,q}|^2=1$, we have $\EE X_{ w(i,j)}=1$ for all $(i,j)\in \Lambda^{(N)}$ such that $[ w(i,j)]=w$.

The above observations yield
\be
\label{eq:momentMeanComb}
\lim_{N\to\infty} \EE m_{k,N}
= \lim_{N\to\infty} \frac{1}{(N\ell)^{k+1}} C^{\lambda}_{k,N},
\ee
where $C^{\lambda}_{k,N}$ is the number of pairs $(T,\mathfrak{c})$, where $T$ is a plane tree on $k+1$ vertices with vertex set $V$ and edge set $E$, and $\mathfrak{c}\colon V\to \llbracket N\ell \rrbracket$ is a labelling of the vertices such that $(\mathfrak{c}(u),\mathfrak{c}(v))\in N\lambda$ for any $\{u,v\}\in E$ and $\mathfrak{c}(u)\neq \mathfrak{c}(v)$ whenever $u,v\in V$ have the same parity (namely, the parity of the graph distance from the root).
Observe now that a function $\mathfrak{c}\colon V\to \llbracket N\ell \rrbracket$ satisfies the property that $(\mathfrak{c}(u),\mathfrak{c}(v))\in N\lambda$ for all $\{u,v\} \in E$ if and only if it is of the form $\mathfrak{c}(v)= (N-1) c(v) + d(v)$ for some labelling $c\colon V\to \llbracket \ell\rrbracket$ such that $(c(u),c(v))\in \lambda$ for any $\{u,v\}\in E$ and for some function $d\colon V\to \llbracket N\rrbracket$.
Notice that one such labelling $c$ makes $(T,c)$ a $\lambda$-plane tree.
Choosing $d$ to be injective (which is possible as long as $N\geq k+1$) or $d$ to be arbitrary, respectively, we obtain the bounds
\[
N(N-1)\cdots (N-k) C^{\lambda}_{k} \leq C^{\lambda}_{k,N}\leq N^{k+1} C^{\lambda}_{k},
\]
where $C^{\lambda}_{k}$ is the number of $\lambda$-plane trees on $k+1$ vertices.
These bounds imply that $C^{\lambda}_{k,N}\sim N^{k+1} C^{\lambda}_{k}$ as $N\to\infty$.
We deduce from~\eqref{eq:momentMeanComb} that $\EE m_{k,N} \to C^{\lambda}_{k}/{\ell^{k+1}} = m_k$ as $N\to\infty$, thus concluding the proof of~\eqref{eq:lem1}.
  
To prove~\eqref{eq:lem2}, we use~\eqref{eq:moment2} to estimate
\[
 \begin{split}
   \Var m_{k,N}
&\leq\frac{1}{(N\ell)^{2k+2}} \sum_{(i,j),(i',j')\in \Lambda^{(N)}} \left|\EE \left[X_{ w(i,j)}X_{ w(i',j')}\right]-\EE X_{ w(i,j)}\EE X_{ w(i',j')}\right| \\
&\leq
\frac{1}{(N\ell)^{2k+2}} \sum_{i,j,i',j'\colon\llbracket k\rrbracket\to\llbracket N\ell\rrbracket} \left|\EE \left[X_{ w(i,j)}X_{ w(i',j')}\right]-\EE X_{ w(i,j)}\EE X_{ w(i',j')}\right|,
\end{split}
\]
since the latter sum contains more (nonnegative) terms than the former.
The last line can be  bounded from above by $O(N^{-2})$, as detailed e.g.\ in~\cite[p.~50]{Bai10}.
This proves~\eqref{eq:lem2}.
  \end{proof}
  
\begin{proof}[Proof of Theorem~\ref{thm:summationFormula}]
The idea behind formula~\eqref{eq:summationFormula} is that the condition~\eqref{eq:lambdaPlaneTrees_equiv} on the labelling of a $\lambda$-plane tree depends on the block-map-values of the vertex labels, rather than on the vertex labels themselves.

Let $\lambda=((a_1+a_2+\cdots+a_r)^{a_1},(a_1+a_2+\cdots+a_{r-1})^{a_2},\ldots,(a_1)^{a_r})$ have $r$ rectangular blocks.
Define the map
\[
\phi_\lambda\colon \{\lambda\text{-plane trees}\} \to \{r\text{-plane trees}\}, \qquad
\phi_\lambda(T,c) := (T',c')
\]
by setting $T':=T$ and $c'(v):=f_\lambda(c(v))$ for every vertex $v$ of $T'=T$.
The map $\phi_\lambda$ is clearly surjective but, unless $\lambda=(r,r-1,\dots,1)$ is the staircase partition, not injective.
More precisely, given an $r$-plane tree $(T',c')$ with $\ell_i$ vertices labelled $i$ for all $i\in\llbracket r\rrbracket$, the preimage $\phi_\lambda^{-1}(T',c')$ has cardinality $a_1^{\ell_1}\cdots a_r^{\ell_r}$, since $f_{\lambda}^{-1}(i)$ has cardinality $a_i$ for all $i\in \llbracket \ell\rrbracket$.
It then follows that the number of $\lambda$-plane trees on $k+1$ vertices is
  \[
  C_{k}^{\lambda}=\sum_{\substack{\ell_1,\ldots,\ell_r\geq 0,\\ \ell_1+\cdots+\ell_r=k+1}} t(\ell_1,\ldots,\ell_r) \, a_1^{\ell_1}\cdots a_r^{\ell_r},
\]
where $t(\ell_1,\dots,\ell_r)$ is the number of $r$-plane trees on $k+1$ vertices with exactly $\ell_i$ vertices labelled $i$ for all $i\in\llbracket r\rrbracket$.
We now use a refined enumeration formula for $r$-plane trees, recently discovered in~\cite[Theorem~1.1]{Okoth24}:
  \be
  \label{eq:OkothWagner}
t(\ell_1,\dots,\ell_r)
=\frac{1}{k}\prod_{j=1}^{\lceil{r/2}\rceil}\binom{\ell_{\geq j}+\ell_{\leq r-j+1}-2}{\ell_j}\prod_{j=1}^{\lfloor{r/2}\rfloor}\binom{{\ell_{\leq j}}+\ell_{\geq r-j+1}-1}{\ell_{r-j+1}},
   \ee
  where $\ell_{\leq i}:=\ell_1+\cdots+\ell_i$ and $\ell_{\geq i}:=\ell_{i}+\cdots+\ell_r$ for any $i\in\llbracket r\rrbracket$.
  It is easy to see that the quantity $t(\ell_1,\dots,\ell_r)$ defined in~\eqref{eq:binomProd} coincides with~\eqref{eq:OkothWagner}.
\end{proof}

\section{$\lambda$-plane trees and $\lambda$-Dyck paths}
\label{sec:lambdaDyckPaths}

In Section~\ref{sec:mainResult} we saw that the limiting spectral distributions of $\lambda$-RMs are characterised by their moments, which are in turn enumerated by $\lambda$-plane trees.
Here we construct a combinatorial bijection between $\lambda$-plane trees and what we call \emph{$\lambda$-Dyck paths}.
These are lattice paths in $\Z^3$, where the first two coordinates move within $\lambda$ (viewed as a subset of $\Z^2$), while the third coordinate performs a usual Dyck path.

It is well known that plane trees on $k+1$ vertices are in bijection with Dyck paths of length $2k$, defined as follows.

\begin{definition}
A \textbf{Dyck path} of length $2k$ is an integer sequence $(h_t)_{t=0}^{2k}$ such that
\begin{enumerate}
\item $h_0=h_{2k}=0$;
\item $h_t-h_{t-1}=\pm 1$ for all $t\in \llbracket 2k\rrbracket$;
\item $h_t\geq 0$ for all $t\in \llbracket 2k-1\rrbracket$.
\end{enumerate}
\end{definition}

\begin{figure}
\begin{subfigure}[c]{0.5\textwidth}
\centering 
\begin{tikzpicture}[scale=0.39]
\draw[very thick,color=Red] (0,0) -- (1,1);
\draw[very thick,color=Blue] (1,1) -- (2,2);
\draw[very thick,color=Blue] (2,2) -- (3,1);
\draw[very thick,color=Red] (3,1) -- (4,0);
\draw[very thick,color=Green] (4,0) -- (5,1);
\draw[very thick,color=Lavender] (5,1) -- (6,2);
\draw[very thick,color=Purple] (6,2) -- (7,3);
\draw[very thick,color=Purple] (7,3) -- (8,2);
\draw[very thick,color=Turquoise] (8,2) -- (9,3);
\draw[very thick,color=Turquoise] (9,3) -- (10,2);
\draw[very thick,color=Lavender] (10,2) -- (11,1);
\draw[very thick,color=gray] (11,1) -- (12,2);
\draw[very thick,color=gray] (12,2) -- (13,1);
\draw[very thick,color=Green] (13,1) -- (14,0);
\draw[very thick,color=Orange] (14,0) -- (15,1);
\draw[very thick,color=Orange] (15,1) -- (16,0);
\foreach \x in {0,...,16}{
	\node at (\x,-0.5) {$\scriptstyle\x$};
	}
\end{tikzpicture}
\end{subfigure}%
\begin{subfigure}[c]{0.5\textwidth}
\centering
\ytableausetup{mathmode, boxsize=2.8em}
\begin{ytableau}
\scriptstyle
{\color{Green}5},{\color{Green}14} & & \scriptstyle {\color{Red}1},{\color{Red}4} & & &\scriptstyle 0,{\color{Orange}15},{\color{Orange}16} \\
\scriptstyle {\color{Lavender}6},{\color{Lavender}11} &\scriptstyle {\color{Turquoise}9},{\color{Turquoise}10} &\scriptstyle {\color{Blue}2},{\color{Blue}3},{\color{Purple}7},{\color{Purple}8} \\
& & \\
\scriptstyle {\color{gray}12},{\color{gray}13} \\
\\
\\
\end{ytableau}
\end{subfigure}%
\caption{An illustration of a $\lambda$-Dyck path $\pi=(i_t,j_t,h_t)_{t=0}^{16}$ of length $16$, with $\lambda=\lambda'=(6^1,3^2,1^3)$ (same Young diagram as in Fig.~\ref{fig:lambda-plane-tree}).
The height of $\pi$ is depicted on the left-hand side: this is a Dyck path $(h_t)_{t=0}^{16}$, where the first and last step of every height excursion are in the same colour.
The location of $\pi$ is depicted on the right-hand side: every box $(i,j)\in\lambda$ contains all `times' $ 0\leq t\leq 16$ such that $(i_t,j_t)=(i,j)$.
The path starts and ends at the same box $(1,6)$.
At even times, the path moves horizontally; at odd times, the path moves vertically.
All the height excursions of $\pi$ are balanced, since the times when a height excursion starts and ends are contained in the same box of $\lambda$ (e.g.\ the red `times' $1$ and $4$ are both in box $(1,3)$).
The corresponding $\lambda$-plane tree $(T,c)$ is the one of Figure~\ref{fig:lambda-plane-tree}, where the edges of $T$ correspond to height excursions in $\pi$ and have the same colour coding.
}\label{fig:lambda-Dyck-paths}
\end{figure}
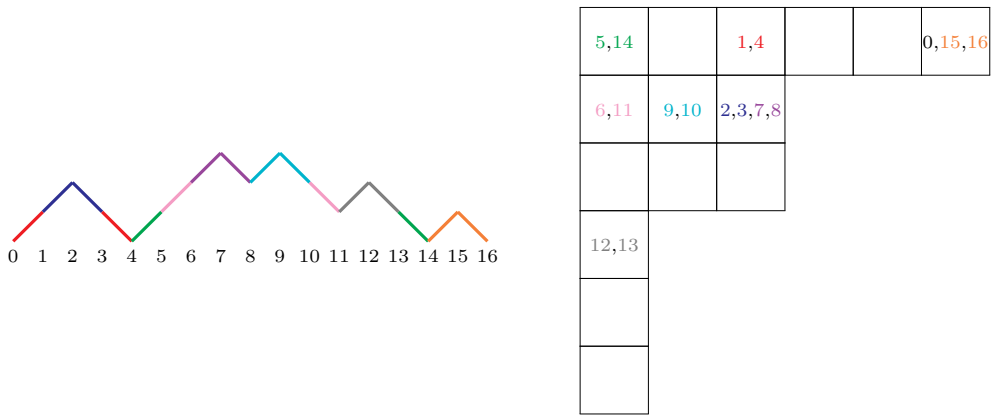

A Dyck path of length $2k$ is usually visualised through the lattice path $(t,h_t)_{t=0}^{2k}$ from $(0,0)$ to $(2k,0)$, completely contained in the quadrant $\{0,1,2,\dots\}\times \{0,1,2,\dots\}$ and consisting of $k$ up-steps (with increment $(1,1)$) and $k$ down-steps (with increment $(1,-1)$).
See Figure~\ref{fig:lambda-Dyck-paths}, left-hand panel.

\begin{definition}
\label{def:lambda-path}
Let $\lambda=\lambda'\in\Par$.
A \textbf{$\lambda$-path} of length $M\geq 0$ is a sequence $\pi=(i_t,j_t,h_t)_{t=0}^M$ in $\Z^3$ such that
\begin{enumerate}
\item\label{def:lambda-path_lambda} $(i_t,j_t)\in\lambda$ for all $t\in\{ 0,\dots,M \}$;
\item\label{def:lambda-path_odd} $i_t-i_{t-1}=0$ for all \emph{odd} $t\in\llbracket M \rrbracket$;
\item\label{def:lambda-path_even} $j_t-j_{t-1}=0$ for all \emph{even} $t\in\llbracket M \rrbracket$;
\item\label{def:lambda-path_height} $h_t-h_{t-1}=\pm 1$ for all $t\in\llbracket M \rrbracket$.
\end{enumerate}
We will say that $(i_t,j_t)$ and $h_t$ are the \textbf{location} and the \textbf{height} of the path at time $t$, respectively.
\end{definition}

In words, the location of a $\lambda$-path moves from a box $(i,j)$ of $\lambda$ to another, alternating between horizontal steps, in which the $i$-coordinate (row of $\lambda$) remains constant, and vertical steps, in which the $j$-coordinate (column of $\lambda$) remains constant; these steps can be of any integer size, including $0$, as long as the location of the path does not exit $\lambda$.
On the other hand, the height of a $\lambda$-path simply moves upwards or downwards by one unit step.

\begin{definition}
A \textbf{height excursion} of a $\lambda$-path $(i_t,j_t,h_t)_{t=0}^M$ is a subpath $(i_t,j_t,h_t)_{t=a}^b$, $0\leq a\leq b\leq M$, such that $h_a=h_b<h_t$ for all $a<t<b$.
A height excursion is said to be \textbf{balanced} if $(i_{a+1},j_{a+1})=(i_b,j_b)$ and \textbf{unbalanced} otherwise.
\end{definition}

In words, a height excursion starts and ends at the same height, with a strictly higher height in between the starting and ending points.
In a balanced height excursion, the first and the last step point towards the same location within $\lambda$.

Notice that, since the height of a $\lambda$-path varies by $\pm 1$ at each time step, the length $b-a$ of a height excursion is necessarily even.

\begin{definition}
A \textbf{$\lambda$-Dyck path} of length $2k$ is a $\lambda$-path $\pi=(i_t,j_t,h_t)_{t=0}^{2k}$ such that
\begin{enumerate}
\item the height $(h_t)_{t=0}^{2k}$ of $\pi$ is a Dyck path;
\item $(i_0,j_0)=(i_{2k},j_{2k})$, i.e.\ the path starts and ends at the same location;
\item $\pi$ contains no unbalanced height excursions.
\end{enumerate}
\end{definition}
See Figure~\ref{fig:lambda-Dyck-paths} for an illustration.

\begin{example}
Consider $\lambda=(1)$, the partition of size $1$.
In such a case, the location of any $\lambda$-Dyck path always stays put at $(1,1)$, i.e.\ any step of the path points towards location $(1,1)$.
In particular, any height excursion is trivially balanced.
Such a $\lambda$-Dyck path is, thus, only determined by its height, which is a Dyck path.
Consequently, when $\lambda=(1)$, $\lambda$-Dyck paths of length $2k$ are in bijection with plane trees on $k+1$ vertices.
\end{example}

The following result extends the simple situation of the latter example to general self-conjugate partitions $\lambda$.

\begin{theorem}
\label{thm:lambdaDyck}
Let $\lambda=\lambda'\in\Par$.
There is a bijection between $\lambda$-plane trees on $k+1$ vertices and $\lambda$-Dyck paths of length $2k$.
\end{theorem}

\begin{proof}
Any (rooted) plane tree $T$ on $k+1$ vertices can be uniquely described via a \textbf{contour walk} on the tree, which determines its root and its planar orientation.
Using the terminology from the proof of Theorem~\ref{thm:main}, a contour walk is an (alternating) spanning closed walk of length $2k$ on the tree.
This is the sequence $ w=(w_0,\ldots,w_{2k})\in V^{2k+1}$, where $V$ is the vertex set of $T$, constructed by starting from the root and following the `contour' of the tree, as follows.
Let $w_0$ be the root.
Inductively, suppose the sequence has been constructed up to $w_{t-1}$.
If $w_{t-1}$ is a leaf, or if all the children of $w_{t-1}$ have been already visited (i.e., for each child $v$ of $w_{t-1}$, there exists $s<t-1$ such that $w_s=v$), then let $w_t$ be the unique parent of $w_{t-1}$.
Otherwise, let $w_t$ be the leftmost child of $w_{t-1}$ that has not been visited yet.
Notice that $\{w_{t-1},w_t\}$ is an edge of $T$ for all $t$.
Let
\[
\epsilon_t:=
\begin{cases}
+1 &\text{if $w_t$ is a child of $w_{t-1}$,} \\
-1 &\text{if $w_{t-1}$ is a child of $w_t$}
\end{cases}
\]
for all $t\in\llbracket 2k\rrbracket$.
In words, $\epsilon_t=+1$ if the walk step from $w_{t-1}$ to $w_t$ is away from the root, and $\epsilon_t=-1$ if it is towards the root.
The sequence $(h_t)_{t=0}^{2k}$ defined by
\begin{equation}
\label{eq:height}
h_0:=0, \qquad h_t:=\epsilon_1+\dots+\epsilon_t \quad\text{for } t\in\llbracket 2k\rrbracket
\end{equation}
also determines uniquely the plane tree $T$ and is the Dyck path of length $2k$ classically associated with $T$.
Notice that $h_t$ is the distance from the root of the vertex $w_t$.

Let now $\lambda=\lambda'\in \Par$ with $\ell:=\ell(\lambda)$.
Take any labelling function $c\colon V \to \llbracket\ell\rrbracket$.
Let $c_t:=c(w_t)$ for all $t\in\{ 0,\dots,2k\}$.
By definition, $(T,c)$ is a $\lambda$-plane tree if and only if $(c_{t-1}, c_t)\in \lambda$ for all $t\in\llbracket 2k\rrbracket$.
Let us now define a sequence $(i_t,j_t)_{t=0}^{2k}$ as follows.
Let $(i_0,j_0):=(c_0,c_{2k-1})$ and, for $t\in\llbracket 2k\rrbracket$, let
\[
(i_t,j_t):=
\begin{cases}
(c_{t-1},c_t) &\text{if $t$ is odd,} \\
(c_t,c_{t-1}) &\text{if $t$ is even.}
\end{cases}
\]

Let us first check that the sequence $\pi=(i_t,j_t,h_t)_{t=0}^{2k}$ is a $\lambda$-path of length $2k$.
For odd $t\in\llbracket 2k\rrbracket$, we have $(i_t,j_t)=(c_{t-1}, c_t)\in \lambda$ (as $(T,c)$ is a $\lambda$-plane tree) and $i_t-i_{t-1}=c_{t-1}-c_{t-1}=0$.
For even $t\in\llbracket 2k\rrbracket$, we have $(i_t,j_t)=(c_t,c_{t-1})\in \lambda$ (as $(c_{t-1}, c_t)\in\lambda$ and $\lambda$ is self-conjugate) and $j_t-j_{t-1}=c_{t-1}-c_{t-1}=0$.
Notice that $w_0=w_{2k}$ is the root, so we have $c_0=c_{2k}$ and
\begin{equation}
\label{eq:initialFinalPoint}
(i_0,j_0)
=(c_0,c_{2k-1})
=(c_{2k},c_{2k-1})
=(i_{2k},j_{2k}),
\end{equation}
which implies $(i_0,j_0)\in\lambda$.
Finally, the fact that $h_t-h_{t-1}=\pm 1$ for all $t\in\llbracket 2k\rrbracket$ follows from~\eqref{eq:height}.
Therefore, all the conditions of Definition~\ref{def:lambda-path} are satisfied with $M=2k$.

We now show that $\pi$ is a $\lambda$-Dyck path of length $2k$.
As already noted above, the height $(h_t)_{t=0}^{2k}$ of $\pi$ is a Dyck path of length $2k$.
The fact that $(i_0,j_0)=(i_{2k},j_{2k})$ follows from~\eqref{eq:initialFinalPoint}.
It remains to check that $\pi$ contains no unbalanced height excursions.
Let $(i_t,j_t,h_t)_{t=a}^b$, $0\leq a\leq b\leq 2k$, be any height excursion, so that $h_a=h_b<h_t$ for all $a<t<b$.
This means that both vertices $w_a$ and $w_b$ are at distance $h_a$ from the root of $T$, while all vertices $w_{a+1},\dots,w_{b-1}$ are at higher distance from the root.
It follows from the construction of the sequence $ w$ that $w_a=w_b$ and $w_{a+1}=w_{b-1}$.
Notice that $a$ and $b$ need to have the same parity, as $h_a=h_b$.
Suppose for example that $a$ is even (the proof when $a$ is odd is similar), so that $b$ is also even.
Then, we have
\[
(i_{a+1},j_{a+1})
=(c_a,c_{a+1})
=(c(w_a),c(w_{a+1}))
=(c(w_b),c(w_{b-1}))
=(c_b,c_{b-1})
=(i_b,j_b).
\]
This proves that the height excursion is balanced.

Let us now briefly explain how the inverse map from $\lambda$-Dyck paths of length $2k$ to $\lambda$-plane trees works.
Let $\pi=(i_t,j_t,h_t)_{t=0}^{2k}$ be a $\lambda$-Dyck path.
Let $T$ be the plane tree associated with the Dyck path $(h_t)_{t=0}^{2k}$ and let $ w$ be the corresponding contour walk.
We need to define the labelling function $c$.
Every vertex $v$ of $T$ is visited at least once by the sequence, so there exists $t$ such that $w_t=v$.
If $t$ is odd, let $c(v):=j_t$, while if $t$ is even let $c(v):=i_t$.
The absence of unbalanced height excursions guarantees that this labelling function $c$ is well defined, in the sense that $c(v)$ is defined in the same way no matter the choice of $t$ such that $w_t=v$.
Moreover, if $\{u,v\}$ is an edge of $T$, then there exists $t\in\llbracket 2k\rrbracket$ such that $u=w_{t-1}$ and $v=w_t$ and, thus, using~\ref{def:lambda-path_odd} and~\ref{def:lambda-path_even} of Definition~\ref{def:lambda-path}, we obtain
\[
(c(u),c(v))=
\begin{cases}
(i_{t-1},j_t)=(i_t,j_t) &\text{if $t$ is odd,} \\
(j_{t-1},i_t)=(j_t,i_t) &\text{if $t$ is even.}
\end{cases}
\]
By~\ref{def:lambda-path_lambda} of Definition~\ref{def:lambda-path} and the fact that $\lambda$ is self-conjugate, we have $(c(u),c(v))\in\lambda$ in both cases.
We conclude that $(T,c)$ is a $\lambda$-plane tree.
\end{proof}

\begin{remark}
Besides the one that we have presented here, other path constructions enumerating $\lambda$-plane trees may exist.
For instance, Gu, Prodinger and Wagner~\cite[Sec. 4]{Gu10} found an explicit `glove bijection' between $r$-plane trees and certain lattice paths, which are a sort of Dyck paths allowing up-steps of size greater than $1$.
It would be interesting to extend such a bijection to $\lambda$-plane trees.
\end{remark}

\section{Limiting spectral distributions: an analytic approach}
\label{sec:generatingFunctions}
The sequence of moments $C^{\lambda}_k/{\ell^{k+1}}$ uniquely identifies the distribution $F^{\lambda}$ of Theorem~\ref{thm:main}, but provides a rather indirect description of it.
The aim of this section is to set up a framework for $\lambda$-RMs to answer rather natural questions in random matrix theory such as: `is $F^{\lambda}$ absolutely continuous?', `what is {its} support?', etc.
More specifically, we work with the generating function of the integer sequence $C^{\lambda}_k$ and prove that various transforms of $F^{\lambda}$ satisfy polynomial equations.
We give explicit solutions in the case of fat hooks, thereby providing a free-probabilistic interpretation of $F^{\lambda}$.

Throughout this section, we fix
\[
\lambda=((a_1+a_2+\cdots+a_r)^{a_1},(a_1+a_2+\cdots+a_{r-1})^{a_2},\ldots,(a_1)^{a_r})\in\Par,
\]
with $\ell:=\ell(\lambda)=a_1+\dots+a_r$.

A first useful observation is that a $\lambda$-RM whose nonzero i.i.d.\ entries have an absolutely continuous distribution is invertible almost surely if and only if the shape $\lambda$ contains the staircase partition $(\ell, \ell-1,\dots,1)$.
More in general, we can compute the dimension of the null space of the random matrix as follows.
 \begin{lemma}
 \label{lem:rank}
 If $X$ is a $\lambda$-RM with nonzero entries that are i.i.d.\ complex continuous random variables, then the dimension of the null space $\operatorname{Ker}X$ is
\begin{align}
\label{eq:dimKer1}
\operatorname{dim}\operatorname{Ker} X
&=\#\left\{i\in\llbracket \ell \rrbracket\colon \lambda_i<\ell-i+1\right\} \\
\label{eq:dimKer2}
&= \sum_{j=2}^r \max\{0,\min\{a_{\geq j} - a_{\leq r-j+1},a_j\}\}
\end{align}
almost surely, where $a_{\leq n}:=a_1+\dots +a_n$ and $a_{\geq n}:=a_n+\dots +a_r$ for any $n\in\llbracket r\rrbracket$.
 \end{lemma}
 \begin{proof}
The proof is based on the fact that, as the nonzero entries are independent and continuous, any linear relation between them has probability zero.
 
 To prove~\eqref{eq:dimKer1}, let $Y$ be the matrix whose entry $(i,j)$ is the  entry $(i,\ell-j+1)$ of $X$.
 Of course we have $\operatorname{rank} Y=\operatorname{rank} X$.
 Perform Gaussian elimination to bring $Y$ in its reduced row echelon form (RREF).
 If $\lambda_i\geq \ell-i+1$ for all $i$, then, almost surely, the RREF of $Y$ is upper triangular with ones on the diagonal, so $Y$ is invertible.
 If $\lambda_i< \ell-i+1$ for some $i$, then, with probability $1$, the RREF of $Y$ has a zero $(i,i)$-entry on the diagonal.
The number of zeros on the diagonal of the RREF of $Y$ is precisely $\operatorname{dim}\operatorname{Ker} Y=\operatorname{dim}\operatorname{Ker} X$.
  
To prove~\eqref{eq:dimKer2}, notice first that, if $f_\lambda$ is the block map associated with the self-conjugate partition $\lambda$, then $f_\lambda(i)=j\in \llbracket r\rrbracket$ if and only if $\lambda_i=a_{\leq r-j+1}$ and $i=a_{\leq j-1}+n$ for some $n\in\llbracket a_j\rrbracket$.
We thus have
\[
\begin{split}
\#\left\{i\in\llbracket \ell \rrbracket\colon \lambda_i<\ell-i+1\right\}
&= \sum_{j=1}^r \#\{n\in \llbracket a_j\rrbracket\colon a_{\leq r-j+1} < \ell -(a_{\leq j-1}+n)+1 \} \\
&= \sum_{j=1}^r \#\{n\in \llbracket a_j\rrbracket\colon n\leq a_{\geq j}-a_{\leq r-j+1} \},
 \end{split}
 \]
using the fact that $a_1+\cdots+ a_r =\ell$.
Noting that the $j=1$ summand is always zero and rewriting
\[
\#\{n\in \llbracket a_j\rrbracket\colon n\leq a_{\geq j}-a_{\leq r-j+1} \}
= \max\{0,\min\{a_{\geq j} - a_{\leq r-j+1},a_j\}\},
\]
we arrive at~\eqref{eq:dimKer2}.
 \end{proof}
 
 The above matrix-theoretic argument implies that $F^{\lambda}(x)$ has a jump (corresponding to an atomic component) at $x=0$ of height at least
 \[
 \frac{1}{\ell}\#\left\{i\in\llbracket\ell\rrbracket \colon \lambda_i<\ell-i+1\right\}.
 \]
 
 We now proceed to study the limiting distribution using some analytic tools.
 We mostly follow the nomenclature of~\cite{Mingo17}, to which we refer the reader for
more background.
Let $F$ be a probability distribution function.
The \textbf{Cauchy transform} of $F$ is
  \[
  G(z):=\int\frac{1}{z-x}\,\de F(x),\qquad \operatorname{Re}z\neq0.
  \]
The distribution $F$ can be recovered from $G$ through the Stieltjes inversion theorem:
  \be
  \label{eq:Stiletjes_inversion}
  F(b)-F(a)=-\lim_{\epsilon\downarrow0}\frac{1}{\pi}\int_{a}^b\operatorname{Im}G^{}(x+i \epsilon)\,\de x,
  \ee
  for $a<b$ such that $a$ is a continuity point of $F^{}$.
  If $F_1$ and $F_2$ are distribution functions with Cauchy transforms $G_1=G_2$, then $F_1=F_2$.
   
The {measure $\de F$} has compact support if and only if its Cauchy transform $G^{}(z)$ has an analytic extension in a neighbourhood of infinity.
Then, denoting by $m_k$, $k\geq0$, the moments of $F$,  
  \[
    G(z)=\sum_{k\geq0}\frac{1}{z^{k+1}}m_k
  \]
  for $|z|$ large enough.
  
{As in~\cite[Sec.~2.3]{Edelman08}, we define the \textbf{$R$-transform} $R(z)$ of $F$ to be} the right inverse, by composition, of  $G$ in a neighbourhood of infinity, i.e.\
\be
\label{eq:def_R}
G(R(z)+1/z)=z,\qquad R(z)=\sum_{k\geq0} R_{k+1} z^{k}.
\ee
The coefficients $R_{k+1}$, $k\geq 0$, are the so-called \textbf{free cumulants} of the distribution $F$.

Given a {measure $\de F$} whose support is contained in $[0,\infty)$, its \textbf{$S$-transform} $S(z)$ is defined by the equation
\be
\label{eq:def_S}
S(z) R(z S(z))=1
\ee
in a neighbourhood of $0$ {(this is equivalent to the definition given in~\cite[Sec.~2.4]{Edelman08}).}

Let  $F_1$ and $F_2$  be probability distribution functions with  $R$-transforms and $S$-transforms given by $R_1$, $S_1$ and $R_2$, $S_2$, respectively.
The \textbf{additive free convolution} of $F_1$ and $F_2$ is the probability distribution $F_1\boxplus F_2$ with $R$-transform $R(z)$ given by
\[
R(z)=R_1(z)+R_2(z).
\]
The \textbf{multiplicative free convolution} of $F_1$ and $F_2$ is defined as the probability distribution $F_1 \boxtimes F_2$ with $S$-transform $S(z)$ given by
\[
S(z)=S_1(z)S_2(z).
\]

Denote now by $M^{\lambda}$ the generating function of the sequence $\left(C_{k}^{\lambda}\right)_{k\geq 0}$, i.e.\
\[
M^{\lambda}(z) := \sum_{k\geq 0} C_{k}^{\lambda} z^{k+1}=\ell z+|\lambda| z^2+C_{2}^{\lambda} z^3+\cdots.
\]
The series converges in a neighbourhood of zero by the comparison test, since $C_{k}^{\lambda}\leq \ell^{k+1}C_k$ (see~\eqref{eq:lambdaPlaneNoBound}) and the generating function of Catalan numbers has positive radius of convergence.
The radius of convergence of $M^{\lambda}$ is strictly less than $1$, since $C^{\lambda}_k\geq 1$ for all $k\geq 0$.
It follows that the Cauchy transform $G^{\lambda}$ of $F^{\lambda}$ has an analytic continuation in a neighbourhood of infinity, and is related to $M^{\lambda}$ via
\be
\label{eq:G_M}
G^{\lambda}(z)={M^{\lambda}\left(\frac{1}{\ell z}\right)}
=\sum_{k\geq0}\frac{1}{z^{k+1}}\frac{1}{\ell^{k+1}}C_{k}^{\lambda}.
\ee
This also shows that {$\de F^{\lambda}$} has compact support.

Let us write $M^{\lambda}=M_1+\dots+M_\ell$, where
\be
\label{eq:genFnDef}
M_i := \sum_{k\geq 0} C_k^{\lambda,i} z^{k+1}, \qquad
C_k^{\lambda,i}:= \#\{\text{$\lambda$-plane trees on $k+1$ vertices with root labelled $i$}\} .
\ee
We can write a system of equations for these generating functions, based on the classical idea that trees can be enumerated recursively.
Take a $\lambda$-plane tree with root labelled $i$ on $k+1$ vertices, $k\geq 1$.
Removing the edge connecting the root $\rho$ to its leftmost child $v$, one obtains two $\lambda$-plane trees, one rooted at $\rho$ (labelled $i$), and one rooted at $v$ (whose label is some $j\leq \lambda_i$).
The total number of vertices of the two trees is, of course, still $k+1$.
It follows that, for $k\geq 1$,
\[
C_k^{\lambda,i} = \sum_{j=1}^{\lambda_i} \sum_{n=0}^{k-1} C_n^{\lambda,i} C_{k-1-n}^{\lambda,j} .
\]
Plugging the latter into~\eqref{eq:genFnDef}, we see that $M_1,\dots,M_\ell$ satisfy the functional equations
\[
M_i = z+M_i\left( M_1 +\dots + M_{\lambda_i}\right), \qquad 1\leq i\leq \ell.
\]

Using multirectangular coordinates for $\lambda=\lambda'$, by~\eqref{eq:selfConjMultirectangular2} we have $\lambda_i=a_1+\dots +a_{r-f_{\lambda}(i)+1}$, where $f_{\lambda}$ is the block map.
If $f_{\lambda}(i)= f_{\lambda}(j)$, i.e.\ $i$ and $j$ are in the same `block' in $\lambda$, then $M_i=M_j$.
Therefore, one may reduce the problem by setting, for all $j\in\llbracket r\rrbracket$, $H_j:= M_i$ for any choice of $i\in\llbracket \ell\rrbracket$ such that $f_{\lambda}(i)=j$.
The reduced collection $H_1,\dots,H_r,M^{\lambda}$ satisfies the polynomial equations
\be
\label{eq:pol_system}
\begin{aligned}
H_1 &= z+H_1 \left(a_1 H_1+\cdots +a_rH_r\right),  \\
H_2 &= z+H_2 \left(a_1 H_1+\cdots +a_{r-1}H_{r-1}\right),  \\
&\;\,\vdots\\
H_r &= z+H_r \left(a_1 H_1\right),  \\
M^{\lambda}&=a_1 H_1+\cdots+a_rH_r.
\end{aligned}
\ee
Solving for $H_1,\dots,H_r$ we obtain the Cauchy transform of $F^{\lambda}$ from~\eqref{eq:G_M} as
\be
\label{eq:Gtrace}
G^{\lambda}(z)={\sum_{j=1}^{r} a_j H_j\left(\frac{1}{\ell z}\right).}
\ee
When $a_1=\dots =a_r=1$, so that $\lambda=(r,r-1,\dots,1)$ is a staircase shape, the polynomial system~\eqref{eq:pol_system} has been solved in~\cite{Gu10}.
For $r=1$, we recover the algebraic equation $M^2-M+z=0$ for the generating function of Catalan numbers.

\begin{remark}
The Cauchy transform of general large dimensional `block' random matrices can be computed using the analytic methods of \emph{operator-valued free probability}.
See~\cite[Ch.~9]{Mingo17} for a general exposition, \cite{Ding14, Arizmendi16} for explicit examples, and~\cite[Sec.~3.3]{Speicher06} for detailed calculations in the case of Wishart type block matrices. 
Adapting the notation {of~\cite{Speicher06}} to our setting, it is possible to construct an $r\times r$ complex-matrix-valued  function $\mathcal{G}^{\lambda}(z)=(g^{\lambda}_{jk}(z))$, analytic for $\operatorname{Im}(z)>0$, whose \emph{weighted trace} coincides with the Cauchy transform:
\[
G^{\lambda}(z)=\tr_{a}\mathcal{G}^{\lambda}(z):={\sum_{j=1}^{r} a_j g^{\lambda}_{jj}(z).}
\]
This is another route to obtain formula~\eqref{eq:Gtrace}, where $g^{\lambda}_{jj}(z)={H_j(1/(\ell z))}$. 

{We also mention that the problem of studying the limiting spectral distribution of $\lambda$-RMs naturally falls in the framework of `rectangular probability spaces' developed in~\cite{Benaych-Georges09}.
In that notation, a random matrix $X_N$ of shape $N\lambda$ is an element of a $(\rho_1,\ldots,\rho_{\ell})$-rectangular probability space with $\rho_1=\cdots=\rho_{\ell}=1/\ell$ all equal and strictly positive.}
\end{remark}

Following Rao and Edelman~\cite{Edelman08}, we say that a probability distribution function $F$ is \textbf{algebraic} if there exists a nonzero bivariate polynomial $L\in \C[G,z]$ such that the Cauchy transform $G(z)$ of $F$ satisfies the equation $L(G(z),z)=0$.
Random matrices whose limiting spectral distribution function is algebraic are referred to as \textbf{algebraic random matrices}.
This class of random matrices is quite special, as the whole information about the limiting spectral distribution is encoded in the bivariate polynomial $L$. 
It turns out that the $(N\lambda)$-RMs that we considered in this article are algebraic.

\begin{theorem}
\label{thm:Stieltjes}
The Cauchy transform $G^{\lambda}(z)$ of $F^{\lambda}$ satisfies the equation
\[
L^{\lambda}(G,z)=0,
\]
where $L^{\lambda}$ is a bivariate polynomial of degree $r+1$ in $G$ and degree $r$ in $z$. 
\end{theorem}

We prove the above result by eliminating the variables $H_1,\ldots,H_r$ from the system of polynomial equations~\eqref{eq:pol_system}.
We mention that the most systematic way to do so would be to compute a Groebner basis and then apply standard results of the elimination theory~\cite[Ch.~3]{Cox15}.
This is easy to perform with the help of computer algebra systems (see~\cite[Appendix~C]{Cox15}).

\begin{proof}
By~\eqref{eq:G_M}, it suffices to show that there exists a bivariate polynomial $P^{\lambda}(M,z)$, of degree $r+1$ in $M$ and of degree $r$ in $z$, such that $P^{\lambda}(M^{\lambda}(z),z)=0$.
Renaming $m,x_1,\ldots,x_r$ the functions $M^{\lambda},a_1 H_1,\ldots, a_rH_r$ respectively, we see from~\eqref{eq:pol_system} that $m,x_1,\ldots,x_r,z$ satisfy the polynomial equations
\begin{align*}
m&=x_1+\cdots+x_r\\
x_1&=x_1(x_1+\cdots+x_r)+a_1z\\
x_2&=x_2(x_1+\cdots+x_{r-1})+a_2z\\
&\;\,\vdots\\
x_{r-1}&=x_{r-1}(x_1+x_2)+a_{r-1}z\\
x_{r}&=x_{r}x_1+a_{r}z
\end{align*}
We will show that it is possible to eliminate the variables $x_1,\ldots,x_r$, obtaining a polynomial equation in $m$ and $z$ only, of degree $r+1$ in $m$ and of degree $r$ in $z$.
Indeed, using the first equation and denoting $x_{\leq i}:=x_1+\dots+ x_i$ and $x_{\geq i}:=x_i+\dots+ x_r$, we can rewrite the remaining $r$ equations as
\begin{align*}
x_1&=x_1m+a_1z\\
x_2&=x_2(m-x_{\geq r})+a_2z\\
x_3&=x_3(m-x_{\geq r-1})+a_3z\\
&\;\,\vdots\\
x_{\lfloor r/2\rfloor}&=x_{\lfloor r/2\rfloor}\left(m-x_{\geq \lceil r/2\rceil+2}\right)+a_{\lfloor r/2\rfloor}z\\
x_{\lfloor r/2\rfloor+1}&=\begin{cases}
x_{\lfloor r/2\rfloor+1}\left(m-x_{\geq \lceil r/2\rceil+1}\right)+a_{\lfloor r/2\rfloor+1}z &r \text{ odd} \\
x_{\lfloor r/2\rfloor+1}x_{\leq \lceil r/2\rceil}+a_{\lfloor r/2\rfloor+1}z &r \text{ even}
\end{cases}\\
&\;\,\vdots\\
x_{r-2}&=x_{r-2}x_{\leq 3}+a_{r-2}z\\
x_{r-1}&=x_{r-1}x_{\leq 2}+a_{r-1}z\\
x_{r}&=x_{r}x_{\leq 1}+a_{r}z
\end{align*}
A suitable rearrangement allows to eliminate the variables sequentially:
\begin{align*}
x_1&=x_1m+a_1z &&\Rightarrow &x_1&=\frac{a_1z}{1-m}\\
x_{r}&=x_{r}x_{\leq 1}+a_{r}z &&\Rightarrow &x_r&=\frac{a_rz}{1-x_{\leq 1}}\\
x_2&=x_2(m-x_{\geq r})+a_2z &&\Rightarrow &x_2&=\frac{a_2z}{1-m+x_{\geq r}}\\
x_{r-1}&=x_{r-1}x_{\leq 2}+a_{r-1}z   &&\Rightarrow &x_{r-1}&=\frac{a_{r-1}z}{1-x_{\leq 2}}\\
x_3&=x_3(m-x_{\geq r-1})+a_3z &&\Rightarrow &x_3&=\frac{a_3z}{1-m+x_{\geq r-1}}\\
x_{r-2}&=x_{r-2}x_{\leq 3}+a_{r-2}z &&\Rightarrow &x_{r-2}&=\frac{a_{r-2}z}{1-x_{\leq 3}}\\
&\;\,\vdots &&\;\;\vdots &&\;\,\vdots \\
x_{\lfloor r/2\rfloor}&=x_{\lfloor r/2\rfloor}\left(m-x_{\geq \lceil r/2\rceil+2}\right)+a_{\lfloor r/2\rfloor}z &&\Rightarrow &x_{\lfloor r/2\rfloor}&=\frac{a_{\lfloor r/2\rfloor}z}{1-m+x_{\geq \lceil r/2\rceil+2}}
\end{align*}
and, finally, according to the parity of $r$:
\begin{align*}
x_{\lfloor r/2\rfloor+1}&=\begin{cases}
x_{\lfloor r/2\rfloor+1}\left(m-x_{\geq \lceil r/2\rceil+1}\right)+a_{\lfloor r/2\rfloor+1}z &r \text{ odd} \\
x_{\lfloor r/2\rfloor+1}x_{\leq \lceil r/2\rceil}+a_{\lfloor r/2\rfloor+1}z &r \text{ even}
\end{cases} \\
\Rightarrow\quad x_{\lfloor r/2\rfloor+1}&=\begin{cases}
\frac{a_{\lfloor r/2\rfloor+1}z}{1-m+x_{\geq \lceil r/2\rceil+1}} &r \text{ odd} \\
\frac{a_{\lfloor r/2\rfloor+1}z}{1-x_{\leq \lceil r/2\rceil}} &r \text{ even}
\end{cases}
\end{align*}
This proves that each variable $x_1,\ldots,x_r$ can be expressed in terms of $m$ and $z$ only.
Inserting these expressions into
\[
m=x_1+x_r+x_2+x_{r-1}+\cdots+x_{\lfloor r/2\rfloor}+x_{\lfloor r/2\rfloor+1}
\]
yields one polynomial equation in $m$ and $z$.
Counting the powers of $m$ and $z$ concludes the proof.
\end{proof}

Recalling~\eqref{eq:def_R}, we immediately deduce a corollary of Theorem~\ref{thm:Stieltjes}.
\begin{corollary}
\label{cor:R}
The $R$-transform $R^{\lambda}$ of $F^{\lambda}$ satisfies the equation
\[
L^{\lambda}(z,R+z^{-1})=0,
\]
which can be reduced to a polynomial equation of degree $r$ in $R$.
\end{corollary}

We now see how the above results can be applied in the case of fat hooks.

We first need to define the \textbf{Marchenko--Pastur distribution} with shape parameter $\alpha>0$ and scale parameter $\beta>0$: this is the measure
\be
\label{eq:MP}
\mathrm{MP}(\alpha,\beta) := \max\left(1-\frac{1}{\alpha},0\right) \delta_0 + \frac{\sqrt{(\gamma_+-x)(x-\gamma_-)}}{2\pi \alpha \beta x} \bm{1}_{[\gamma_-,\gamma_+]}(x) \,\de x,
\ee
where {$\gamma_{\pm}:= \beta(1\pm\sqrt{\alpha})^2$}, and $\de x$ denotes the Lebesgue measure on $\R$.
It can be shown, using the explicit formulae for the moments of $\mathrm{MP}(\alpha,\beta)$ together with~\eqref{eq:def_R} and~\eqref{eq:def_S}, that the $R$-transform and the $S$-transform of $\mathrm{MP}(\alpha,\beta)$ are
\be\label{eq:R&S_MP}
{
R_{\mathrm{MP}(\alpha,\beta)}(z)= \frac{\beta}{1-\alpha\beta z}, \qquad\quad
S_{\mathrm{MP}(\alpha,\beta)}(z)= \frac{1}{\beta(1+\alpha z)}.}
\ee

We also denote the Bernoulli distribution of parameter $p\in [0,1]$ as
\be\label{eq:Bernoulli}
\mathrm{Ber}(p) := (1-p)\delta_0 + p \delta_1 .
\ee
It can be readily checked from~\eqref{eq:def_R} and~\eqref{eq:def_S} that its $R$-transform and its $S$-transform are
\be\label{eq:R&S_Ber}
R_{\mathrm{Ber}(p)}(z)= \frac{z-1+\sqrt{1-2(1-2p)z+z^2}}{2z}, \qquad\quad
S_{\mathrm{Ber}(p)}(z)= \frac{1+z}{p+z}.
\ee

\begin{theorem}
\label{thm:fatHook}
Let $\lambda=((a_1+a_2)^{a_1},a_1^{a_2})$ be a fat hook of length $\ell:=a_1+a_2$.
{Let $p_1:=a_1/\ell$ and $p_2:=a_2/\ell$.}
Then:
\begin{enumerate}
\item The Cauchy transform $G^{\lambda}(z)$ of $F^{\lambda}$ satisfies $L^{\lambda}(G^{\lambda},z)=0$, where $L^{\lambda}$ is the bivariate polynomial
\be
\label{eq:G_fathook}
L^{\lambda}(G,z):={z^2 G^3 + z\left(p_1  -p_2  -2 z\right)G^2+ z \left(2p_2+ z\right) G+ \big(p_1^2- z\big).}
\ee
\item The $R$-transform of $F^{\lambda}$ is given by
\be
\label{eq:R_fathook}
R^{\lambda}(z)
= {\frac{p_1}{1- z} +\frac{1}{2z} \left(\sqrt{\frac{1-\left(p_2-p_1\right)^2z}{1-z}} -1\right).}
\ee
\item The measure $\de F^{\lambda}$ can be expressed in terms of free convolutions of Marchenko--Pastur and Bernoulli distributions as follows:
\begin{align}
\label{eq:MP-BerConv1}
\de F^{\lambda}
&= \mathrm{MP}\left({p_1^{-1},p_1}\right) \boxplus \left(\mathrm{MP}\left({p_1^{-1},p_1}\right) \boxtimes \mathrm{Ber}\left({p_2}\right) \right) \\
\label{eq:MP-BerConv2}
&= \mathrm{MP}\left({p_1^{-1},p_1}\right) \boxplus \left(\mathrm{MP}\left({p_2^{-1},p_2}\right) \boxtimes \mathrm{Ber}\left({p_1}\right) \right)
\end{align}
\end{enumerate}
\end{theorem}

In the special case $a_1=a_2=1$, the above results match those obtained in~\cite{Mlotkowski13} with different methods (notice that their definition of the $R$-transform is slightly different from ours).

\begin{proof}
\begin{enumerate}
\item
For $r=2$, the polynomial equations~\eqref{eq:pol_system} for the reduced generating functions read
\[
\begin{aligned}
H_1 &= z+H_1(a_1H_1 + a_2H_2), \\
H_2 &= z+a_1H_2H_1, \\
M^{\lambda} &=a_1 H_1+a_2H_2.
\end{aligned}
\]
Eliminating the variables $H_1$ and $H_2$, as detailed in the proof of Theorem~\ref{thm:Stieltjes}, we find that $M^{\lambda}$ solves
\[
M^3 +  ( (a_1-a_2)z-2)M^2 +  (1 + 2 a_2 z)M - z (\ell - a_1^2 z)=0.
\]
It follows from~\eqref{eq:G_M} that the Cauchy transform $G^{\lambda}(z)$ satisfies $L^{\lambda}(G^{\lambda},z)=0$, where $L^{\lambda}(G,z)$ is the bivariate polynomial defined in~\eqref{eq:G_fathook}.
Notice that, in accordance with Theorem~\ref{thm:Stieltjes}, $L^{\lambda}(G,z)$ is of degree $r+1=3$ in $G$ and $r=2$ in $z$.
\item
By Corollary~\ref{cor:R}, the $R$-transform of $F^{\lambda}$ solves $a(z) R^2+ b(z)R+c(z)=0$, with
\begin{align*}
a(z)&:= {z(1- z)^2,} \\
b(z)&:= {(1- z)(1-(2p_1+1)z),} \\
c(z)&:= {p_1(2 z+p_1-2).}
\end{align*}
Solving the quadratic equation in $R$, we arrive at~\eqref{eq:R_fathook} (the choice of the root is such that $R^{\lambda}(z)$ is analytic at $z=0$; see~\eqref{eq:def_R}).
\item
By~\eqref{eq:R_fathook}, we have $R^{\lambda}=R_1 + R_2$, where
\begin{align*}
R_1(z) := {\frac{p_1}{1- z}}, \qquad
R_2(z) := {\frac{1}{2z} \left(\sqrt{\frac{1-\left(p_2-p_1\right)^2z}{1-z}} -1\right)}.
\end{align*}
By~\eqref{eq:R&S_MP}, the first term is the $R$-transform of $\mu_1:=\mathrm{MP}\left({p_1^{-1},p_1}\right)$.
If $R_2$ is the $R$-transform of a measure $\mu_2$, then the $S$-transform $S_2$ of $\mu_2$ satisfies~\eqref{eq:def_S}.
Solving for $S$ and recalling that $p_1+p_2=1$, we obtain
\[
\begin{split}
S_2(z) &= {\frac{4(1+z)}{(1+2z-p_1+p_2)(1+2z+p_1-p_2)}} \\
&= {\frac{1}{p_1(1+ z/p_1)} \cdot \frac{1+z}{p_2+z}
= S_{\mathrm{MP}(p_1^{-1},p_1)}(z) \cdot S_{\mathrm{Ber}(p_2)}(z) ,}
\end{split}
\]
using formulas~\eqref{eq:R&S_MP} and~\eqref{eq:R&S_Ber} for the $S$-transforms of the Marchenko--Pastur and Bernoulli distributions.
Therefore, by definition of multiplicative free convolution, $R_2$ is the $R$-transform of the measure $\mu_2:=\mathrm{MP}\left({p_1^{-1},p_1}\right) \boxtimes \mathrm{Ber}\left({p_2}\right)$.
By definition of additive free convolution, we have $\de F^{\lambda}= \mu_1\boxplus \mu_2$, i.e.~\eqref{eq:MP-BerConv1}.
As $R_2$ is symmetric in {$p_1$ and $p_2$}, \eqref{eq:MP-BerConv2} also follows.
\qedhere
\end{enumerate}
\end{proof}

\begin{remark}
\label{rem:classicalRMs}
Theorem~\ref{thm:fatHook} provides an appealing alternative description of the limiting spectral measure $F^{\lambda}$ for fat {hooks} in terms of sums and products of `classical' (i.e., \emph{non}-$\lambda$-shaped) random matrix ensembles.
Consider the following sequences of matrices: let $A_N$ and $B_N$ be $(N\ell)\times (Na_1)$ independent random matrices whose entries are independent standard complex Gaussian, and let $D_N$ be the $(N\ell) \times (N\ell)$ diagonal (non-random) matrix
\[
D_N=\operatorname{diag}(\underbrace{0,\ldots,0}_{\text{$a_1$ times}},\underbrace{1,\ldots,1}_{\text{$a_2$ times}},\ldots,\underbrace{0,\ldots,0}_{\text{$a_1$ times}},\underbrace{1,\ldots,1}_{\text{$a_2$ times}}).
\]
The spectral distributions of the sequences $\frac{1}{N{\ell}}A_N A_N^*$ and $\frac{1}{N{\ell}}B_N B_N^*$ both converge to $\mathrm{MP}\left({\frac{\ell}{a_1},\frac{a_1}{\ell}}\right)$, while the spectral distribution of the sequence $D_N$ is $\mathrm{Ber}\left(\frac{a_2}{\ell}\right)$.
These sequences are asymptotically free, hence~\cite{Mingo17} the spectral distribution of 
\[
\frac{1}{N{\ell}}A_N A_N^*+\frac{1}{N{\ell}}D_NB_N (D_N B_N)^*
\]
converges to the measure on the right-hand side of~\eqref{eq:MP-BerConv1}, which equals $\de F^{\lambda}$.
\end{remark}

{
It would be certainly interesting to see whether formulae similar to~\eqref{eq:MP-BerConv1}-\eqref{eq:MP-BerConv2}, involving free convolutions of simpler distributions, hold for dilations of arbitrary partitions $\lambda$.
This could provide alternative descriptions of $F^{\lambda}$ in terms of sums and products of `classical' random matrix ensembles, in the spirit of Remark~\ref{rem:classicalRMs}.
}

We now proceed to find an explicit expression for the limiting distribution $F^{\lambda}$, with $\lambda$ a fat hook as before.
We essentially use the method explained in~\cite{Edelman08}, to which the reader is referred.

First, notice from~\eqref{eq:G_fathook} that $z=0$ is a singularity of $L^{\lambda}$, in the sense that, at $z=0$, the degree of $L^{\lambda}$ as a polynomial in $G$ decreases.
Such a singularity is the hallmark~\cite{Edelman08} of a possibly nonzero atomic component of the limiting distribution at $x=0$.
Without loss of generality, we can assume that the nonzero entries of $X_N$ have an absolutely continuous distribution.
Hence, from Lemma~\ref{lem:rank} we have  
\[
\dim\operatorname{Ker}(X_N)=\max\{{(a_2-a_1)N},0\}\quad\text{a.s.},
\]
so that the spectrum of $W_N$ will have $\max\{a_2-a_1,0\}N$ zero eigenvalues, a.s..
It follows that the limiting measure $\de F^{\lambda}$ has an atomic part at zero of mass {$\max\left\{p_2-p_1,0\right\}$}, {where $p_1:=a_1/\ell$ and $p_2:=a_2/\ell$ are defined as in Theorem~\ref{thm:fatHook}.} 

To find the full support of the measure $\de F^{\lambda}$, by the Stieltjes inversion theorem~\eqref{eq:Stiletjes_inversion},  it is sufficient to consider $z$ real and find the set $\{z\in\R\colon D_{L^{\lambda}}(z)\leq0\}$, where $D_{L^{\lambda}}(z)$ is the discriminant of $L^{\lambda}(G,z)$ considered as a polynomial in $G$.
Indeed, $D_{L^{\lambda}}(z)$ is proportional to the squared differences of the three roots of $L^{\lambda}(\cdot,z)${, considered as a polynomial in $G$}.
{If the three roots are all real} at $z$, then $G^{\lambda}(z)$ is real and, by the Stieltjes inversion theorem, $z$ does not belong to the support of the measure.
Therefore, a necessary condition for $z$ to belong to the support of the measure is that $D_{L^{\lambda}}(z)<0$, so that two roots are complex conjugate with nonzero imaginary part.

Computing the discriminant
\be
\label{eq:discrim}
D_{L^{\lambda}}(z) =
{z^3  p_1^2 \left(4p_2z^2+\left(p_1^2-20 p_1p_2-8 p_2^2\right)z-4 \left(p_1- p_2\right)^3\right),}
\ee
we see that $D_{L^{\lambda}}(z)\leq0$ if and only if either $z=0$ or $z\in [z_-,z_+]$, where
\[
{z_{\pm}=\frac{8 p_2^2+20p_1p_2-p_1^2\pm(p_1+8p_2)\sqrt{p_1^2+8p_1p_2}}{8 p_2}.}
\]
Therefore, $G^{\lambda}(z)$ has the form
\be
G^{\lambda}(z)={\max\left\{{p_2-p_1},0\right\}\frac{1}{z}+\tilde{G}^{\lambda}(z),}
\ee 
where $\tilde{G}^{\lambda}(z)$ is the Cauchy transform of a sub-probability measure {$\de \tilde{F}^{\lambda}$} with support contained in $[z_-,z_+]$, and $\tilde{G}^{\lambda}(z)\sim\left(1-{\max\left\{p_2-p_1,0\right\}}\right)\frac{1}{z}$ as $z\to\infty$.
Recall also that $F^{\lambda}$ (and hence $\tilde{F}^{\lambda}$) is the distribution function of a measure supported on $[0,\infty)$. 

\begin{figure}
	\centering
	\ytableausetup{boxsize=.5em}
	\includegraphics[width=.45\textwidth]{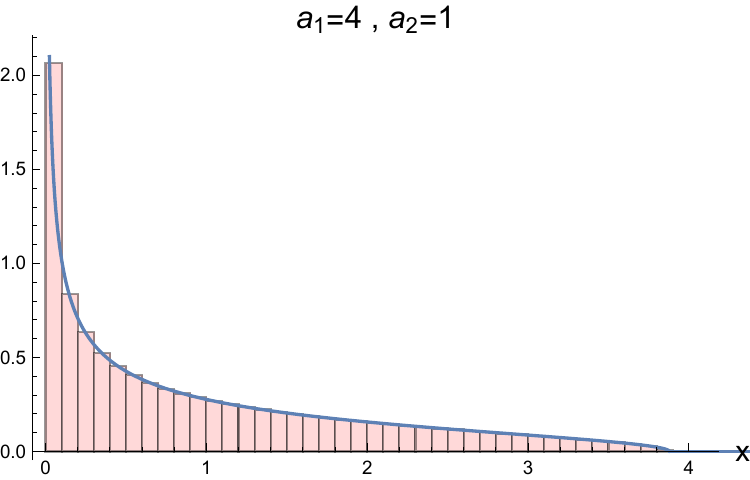}
	\llap{\raisebox{3cm}{ \ydiagram{5,5,5,5,4} \hspace{1.5cm}}}
	\qquad
	\includegraphics[width=.45\textwidth]{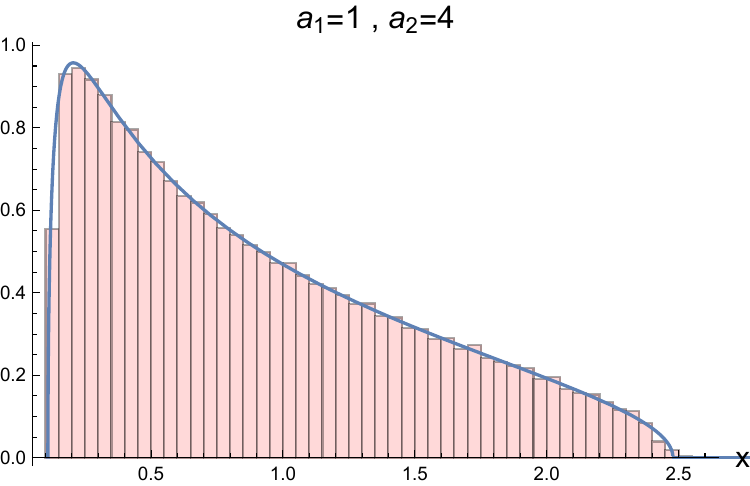}
	\llap{\raisebox{3cm}{ \ydiagram{5,1,1,1,1} \hspace{1.5cm}}}
\caption{Graph of $f^{\lambda}(x)$, the density~\eqref{eq:density_fathook} of the non-atomic part of the limiting spectral measure~\eqref{eq:lim_distr_fat_hook}, compared with numerical data (from a sample of $1000$ matrices with $N=30$), for two fat hooks $\lambda$ of length $\ell(\lambda)=a_1+a_2=5$: the notched square $(5^4,4^1)$ on the left and the hook shape $(5^1,1^4)$ on the right.
The atomic part of the measure for the hook shape at $x=0$ is not shown here.}
\label{fig:varie}
\end{figure}

Applying the Stieltjes inversion theorem to $\tilde{G}^{\lambda}$, we find that the limiting distribution $F^{\lambda}(x)$ corresponds to a sum of a discrete measure at $x=0$ and an absolutely continuous measure:
\be
\label{eq:lim_distr_fat_hook}
\de F^{\lambda}(x)=\max\left\{{p_2-p_1},0\right\}\delta_0+f^{\lambda}(x)\,\de x.
\ee 
The density $f^{\lambda}(x)$ is supported in the interval $[z_-,z_+]\cap[0,+\infty)$ and is given by
\be
\label{eq:density_fathook}
{f^{\lambda}(x)=\frac{1}{\pi  2^{4/3}x^{4/3}}\left| \left| P_{\lambda}(x) +\sqrt{-D_{L^{\lambda}}(x)} \right|^{\frac{1}{3}}\!+\!\frac{2^{2/3}3^{-1} x^{2/3} \left(2 p_1 (p_2+2 x)-p_1^2-(p_2-x)^2\right)}{ \left| P_{\lambda}(x) +\sqrt{-D_{L^{\lambda}}(x)} \right|^{\frac{1}{3}}}\right|,}
\ee
where 
\[
{P_{\lambda}(x)=3^{-3/2} x \left(2 (p_2-x)^3-6 p_1  (p_2-x) (p_2+2 x)+3p_1^2 (2 p_2-5 x)-2 p_1^3\right)}
\]
and $D_{L^{\lambda}}(x)$ is the discriminant~\eqref{eq:discrim}. 
Note that  $z_{-}>0$ if and only if {$p_1< p_2$ (i.e.\ $p_1<1/2$)}, in which case there is a gap between the atomic part at $x=0$ and the support of the density $f^{\lambda}$.
See Fig.~\ref{fig:varie} for a numerical illustration.




\providecommand{\bysame}{\leavevmode\hbox to3em{\hrulefill}\thinspace}
\providecommand{\MR}{\relax\ifhmode\unskip\space\fi MR }
\providecommand{\MRhref}[2]{%
  \href{http://www.ams.org/mathscinet-getitem?mr=#1}{#2}
}
\providecommand{\href}[2]{#2}

\begin{acks}
This project was partially carried out during two research stays of FDC at TU Wien and a research stay of EB at the University of Bari.
FDC is grateful to the Probability Research Unit of TU Wien for hospitality and for providing a stimulating working environment.
EB thanks the Mathematical Physics group at the University of Bari, in particular Marilena Ligab\`o and Giovanni Gramegna, for friendly and interesting discussions.
FDC is grateful to Roberto La Scala for consultations in commutative algebra.
\end{acks}


\end{document}